\documentclass[preprint,10pt]{elsarticle}
\usepackage{cjhebrew}
\usepackage{bm}
\usepackage{xcolor}

\makeatletter
\def\ps@pprintTitle{%
  \let\@oddhead\@empty
  \let\@evenhead\@empty
  \let\@oddfoot\@empty
  \let\@evenfoot\@oddfoot
}
\makeatother

\usepackage{amssymb}

\newtheorem{theorem}{Theorem}[section]
\newtheorem{proposition}{Proposition}[section]
\newtheorem{corollary}{Corollary}[section]

\newtheorem{example}{Example}

\newproof{proof}{Proof}

\begin{document}

\begin{frontmatter}


\title{Modifying the Field Axioms to Create Infinite and Infinitesimal Real Numbers}
\author{Brendan Santangelo} 
\address{brendan.santangelo@gmail.com}

\begin{abstract}
Motivated by the algebraic impossibility of division by zero, a new algebraic structure called an ascended field is introduced by modifying the standard field axioms.  An arithmetic on tuples containing elements from ascended fields is developed, leading to the development of a field extension based around the quotient space of an equivalence on elements in an ascended field.  An extended version of this framework is constructed based on a particular subspace of a quotient space containing tuples of elements from ascended fields.  Structure uniqueness, the ordering of elements, and the totality of the basic algebraic operations are thoroughly explored within these quotient spaces.  This culminates in the development of a unique structure called a complete s-extension of the reals that extends the standard real numbers, maintains totality of the basic arithmetic operations for all elements in the extension, and contains what are defined as infinite and infinitesimal real numbers. 
\end{abstract}

\begin{keyword}
ascended fields \sep modified field axioms \sep division by zero \sep \newline factorials of negative integers \sep infinite real numbers \sep infinitesimal real numbers
\MSC[2020] 08A05, 08A99
\end{keyword}

\end{frontmatter}

\section*{Introduction} 
It is well known that division by zero is an undefined operation.  Interestingly though, the reason this operation is undefined has little to do with division and everything to do with multiplication.  If we could define, say, $\frac{1}{0}$ in a field and set $\frac{1}{0}=x$, then as a consequence of division, we would have $0\cdot x=1$.  This is impossible in a field though as $0\cdot x$ always equals $0$.  But what if multiplication by zero didn't always result in zero and could instead result in $1$, $2$, or even $\pi$?  Would division by zero become possible? 
\newline \indent This provides the motivation to look for a new structure which allows for multiplication by zero to produce results other than zero.  In Section 1, a structure called an ascended field is defined by modifying the standard field axioms to achieve this desired feature.  However, this is not enough to allow for division like $\frac{1}{0}$.  There can now be multiple results for this operation, resulting in an indeterminacy issue.  
\newline \indent To give division by zero a unique result, in Section 3, a certain equivalence relation is defined for elements in ascended fields.   
We call the quotient space of this particular equivalence relation the simple s-extension of a field.  Simple s-extensions uniquely extend any field.  Furthermore, division by zero becomes possible in simple s-extensions and every element that is divided by zero results in its own, uniquely different element from this extension.
\newline \indent If the field of real numbers is extended into its simple s-extension, then the elements in this extension can be ordered.  After they are ordered, we find that all of the elements in this extension which can be the result of division by zero are infinite numbers as each are either greater than or less than every possible real number.  In the simple s-extension of the reals, we also find that factorials for negative integers can be defined.  In fact, we can prove the general formula $(-n)!=[{(-1)^{(n-1)}\over (n-1)!}]$ which gives us a uniquely different infinite output for the factorial of any negative integer $-n\in\Bbb Z$. 
\newline \indent But while division by zero becomes possible and results in infinite numbers, the basic arithmetic operations in a simple s-extension are very limited and are not defined for all elements.  So in Section 2, we introduce tuples whose components are elements from ascended fields.  We consider both tuples containing a finite amount of components as well as tuples containing countably infinite amounts of components.  An arithmetic for these tuples is developed and an equivalence relation for these tuples is defined.  In Section 4, we then consider a subspace of the resulting quotient space called the complete s-extension of a field.  We will find that complete s-extensions uniquely extend any field, act as an extension to non-zero elements in simple s-extensions, and maintain totality of the basic arithmetic operations for all elements in the extension.  However, the zero element, which we could divide by in simple s-extensions, is replaced by a new zero element called absolute zero in complete s-extensions.  
\newline \indent If the field of real numbers is extended into its complete s-extension, then the elements in this extension can be ordered.  After they are ordered, we find that there exist not just finite numbers which correspond to the non-zero reals, but also infinite numbers of different levels as well as infinitesimal numbers of different levels.  Infinite numbers are again greater than or less than every possible real number while infinitesimal numbers are not absolute zero, but are less than every positive real number and greater than every negative real number.  
\newline \indent We call these Infinite and Infinitesimal Real Numbers and its worth noting that each of these numbers is uniquely identified by a non-zero real number value and an integer level.  Any number whose level is $0$ is a finite number, any number whose level is positive is an infinite number, and any number whose level is negative is an infinitesimal number.  So we can think of this arithmetic system as a system in which non-zero, real number values can take on different forms depending on their level.  It makes sense to think of a non-zero value $x\in\Bbb R$ as being in its base form, or finite form, when the level is $0$.  We can think of a non-zero value as being in its first infinite form when the level is $1$ and in its first infinitesimal form when the level is $-1$.
\newpage

\section{Ascended Structures}
Let $(\Bbb S,+,\cdot)$ be an algebraic structure where $+$ is addition and $\cdot$ is multiplication which are well defined and satisfy 
\begin{list}{}{}
\item {1) $(\Bbb S,+)$ forms a commutative group with identity $0$ and inverses $-s$.}
\item {2) $(\Bbb S,\cdot)$ is closed.}
\item {3) $0\cdot 0=0$.}
\item {4) There exists $s\in \Bbb S$ such that $0\cdot s\not=0$ or $s\cdot 0\not =0$ (or both).}
\end{list}
We call $(\Bbb S,+,\cdot)$ an {\bf Ascended Structure}, or just an {\bf S-Structure}.  Since $(\Bbb S,+)$ forms a commutative group, we can define subtraction in an ascended structure $(\Bbb S,+,\cdot)$ with elements $s,t\in\Bbb S$ as
$$s-t=s+(-t).$$
Now if $(\Bbb S,+,\cdot)$ is an ascended structure containing elements $s,t\in\Bbb S$, then we can define an equivalence between $s$ and $t$ by saying that 
$$s=t\,\,\,{\rm if\,\,and\,\,only\,\,if}\,\,\,s-t=0.$$
And since $(\Bbb S,+)$ forms a commutative group, we know that 
$=$ is an equivalence relation.  Finally, we say an ascended structure $(\Bbb S,+,\cdot)$ is {\bf Commutative} if $(\Bbb S,\cdot)$ is commutative.
\begin{example} 
\rm Let $(\Bbb S,+,\cdot)=(\Bbb R\times \Bbb R,+,\cdot)$ where $(\Bbb R,+,\cdot)$ is the reals and
$$\begin{array}{rrcll}
s+t&=&(x,y)+(u,v)&=&(x+u\,,\,y+v). \cr
s\cdot t&=&(x,y)\cdot(u,v)&=&(xu+y+v,\,yv). \cr
\end{array}$$
\end{example}
The structure $(\Bbb R\times \Bbb R,+,\cdot)$ is a commutative S-Structure where $0=(0,0)$ is the zero element satisfying $(0,0)\cdot(x,y)=(y,0)$.  Note that for all $y\not=0$, we have $(0,0)\cdot(x,y)\not=(0,0)$.  
\newline
\newline \indent For every commutative S-Structure $(\Bbb S,+,\cdot)$, we define the subset $\Bbb S_0\subset\Bbb S$ such that 
$$\Bbb S_0=\{x\in\Bbb S\,|\,0\cdot x=0\}.$$
We call an element $x\in\Bbb S_0$ a {\bf Scalar}.  In general, we define the subsets $\Bbb S_{\alpha}\subset\Bbb S$ where $\alpha\in\Bbb S$ is any element in $\Bbb S$ as 
$$\Bbb S_{\alpha}=\{s\in {\Bbb S}\,|\, 0\cdot s=\alpha \}.$$
Furthermore, we define the subset $\Lambda\subseteq{\Bbb S}$, called the {\bf Index of $({\Bbb S},+,\cdot)$}, as
$$\Lambda=\{\alpha\in{\Bbb S}\,|\,{\Bbb S}_{\alpha}\not=\emptyset\}.$$
The third and fourth requirement in the definition of s-structures guarantees that there will be at least two subsets $\Bbb S_\alpha$ which are non-empty as there are at least two elements $\alpha\in\Bbb S$ which satisfy $0\cdot s=\alpha$ for some $s\in\Bbb S$: $\alpha=0$ and some $\alpha\not=0$.  The index $\Lambda$ tells us all of the possible elements $\alpha\in\Bbb S$ such that $0\cdot s=\alpha$ for some $s\in\Bbb S$.  

\begin{proposition}
Let $({\Bbb S},+,\cdot)$ be a commutative S-Structure, then $\{{\Bbb S}_\alpha\}_{\alpha\in\Lambda}$ is a partition on ${\Bbb S}$.
\end{proposition}

\begin{proof}
Let $s_1\in \Bbb{S}$, since $({\Bbb S},\cdot)$ is closed, we know that $0\cdot s_1=\alpha$ for some $\alpha \in \Bbb{S}$.  This means that $s_1\in \Bbb{S}_\alpha$ and since ${\Bbb S}_{\alpha}\not=\emptyset$, we see that $\alpha\in\Lambda$.  Therefore $s_1\in \bigcup_{\alpha\in \Lambda} \Bbb{S}_\alpha$, which means $\Bbb{S} \subseteq \bigcup_{\alpha\in \Lambda} \Bbb{S}_\alpha$.  Furthermore, we see that $\bigcup_{\alpha\in \Lambda} \Bbb{S}_\alpha\not=\emptyset$.  If instead $s_2\in \bigcup_{\alpha\in \Lambda} \Bbb{S}_\alpha$, then $s_2\in \Bbb{S}_\alpha$ for some $\alpha\in \Lambda$.  But ${\Bbb S}_{\alpha}\subset{\Bbb S}$, implying $s_2\in {\Bbb S}$.  Therefore $\bigcup_{\alpha\in \Lambda} \Bbb{S}_\alpha\subseteq {\Bbb S}$, proving that $\Bbb{S}=\bigcup_{\alpha\in \Lambda} \Bbb{S}_\alpha$.
\newline \indent Now let $\alpha_1,\alpha_2\in {\Bbb S}$.  If ${\alpha_1}={\alpha_2}$, then $\Bbb{S}_{\alpha_1}\cap \Bbb{S}_{\alpha_2} =\Bbb{S}_{\alpha_1}\cap \Bbb{S}_{\alpha_1} =\Bbb{S}_{\alpha_1}$.  If ${\alpha_1}\not={\alpha_2}$, then $\Bbb{S}_{\alpha_1}\cap \Bbb{S}_{\alpha_2}=\{s\in {\Bbb S}\,|\,(0 \cdot s)={\alpha_1}={\alpha_2}\}=\emptyset$.  So each ${\Bbb S}_\alpha$ is disjoint from one another.  Therefore $\{{\Bbb S}_\alpha\}_{\alpha\in\Lambda}$ is a partition on ${\Bbb S}$. $\square$
\end{proof} 
For a commutative S-Structure $(\Bbb S,+,\cdot)$, we say the following:
\begin{list}{}{}
\item {1) If $\Lambda\subseteq \Bbb S_0$, we say that $(\Bbb S,+,\cdot)$ is {\bf Complete}.}
\item {2) If $\Bbb S_0\subseteq \Lambda$, we say that $(\Bbb S,+,\cdot)$ is {\bf Regular}.}
\item {3) If $\Lambda=\Bbb S_0$, we say that $(\Bbb S,+,\cdot)$ is {\bf Complete Regular}.}
\end{list}    

It is rather straight-forward to show that every S-Structure lacks the distributive property.

\begin{proposition}
Let $(\Bbb S,+,\cdot)$ be an S-Structure, then $(\Bbb S,+,\cdot)$ is not distributive.
\end{proposition}

\begin{proof}
If we assume $({\Bbb S},+,\cdot)$ is an S-Structure with distribution and $s\in {\Bbb S}$ such that $0\cdot s\not=0$, then $0\cdot s=(0+0)\cdot s=(0\cdot s)+(0\cdot s)$.  By adding $-(0\cdot s)$ to both sides, we get that $0\cdot s=0$.  But we know that $0\cdot s\not=0$, giving us a contradiction.  If instead we have $s\in\Bbb S$ such that $s\cdot 0\not=0$, then we can similarly conclude that $s\cdot 0=0$.  This again gives us a contradiction, proving that $({\Bbb S},+,\cdot)$ is not distributive. $\square$
\end{proof}
Instead, we say that a commutative S-Structure $(\Bbb S,+,\cdot)$ is {\bf Wheel Distributive} if for every $s,t,r\in\Bbb S$, we have 
$$s\cdot (t+r)+(s\cdot 0)=(s\cdot t)+(s\cdot r).\,\,\,\,\,\,  {\rm (See\,reference}\,\cite{car},\,{\rm page\,5\,\,of\,article)}$$  
Equivalently, we can use $s\cdot (t+r)=(s\cdot t)+(s\cdot r)-(s\cdot 0)$.  Furthermore, we say that a commutative S-Structure $(\Bbb S,+,\cdot)$ is {\bf S-Associative} if for every $s\in\Bbb S$ and every $x,y\in\Bbb S_0$, we have 
$$x\cdot (y\cdot s)+\big(\big((x-1)\cdot(y-1)\big)\cdot(0\cdot s)\big)=(x\cdot y)\cdot s.$$  
Equivalently, we can use $x\cdot (y\cdot s)=(x\cdot y)\cdot s-\big(\big((x-1)\cdot(y-1)\big)\cdot(0\cdot s)\big)$.  Finally, we say an S-Structure has {\bf Scalar Inverses} if for every $x\in\Bbb S_0$ where $x\not=0$, there exists $x^{-1}\in\Bbb S_0$ such that $x\cdot x^{-1}=x^{-1}\cdot x=1$.
 
\section*{Ascended Rings and Ascended Fields}
Let $(\Bbb S,+,\cdot)$ be an algebraic structure with well defined addition and multiplication operations.  We say that $(\Bbb S,+,\cdot)$ is an {\bf Ascended Ring} if
\begin{list}{}{}
\item {1) $(\Bbb S,+)$ is a commutative group with identity $0$ and inverses $-s$.}
\item {2) $(\Bbb S,\cdot)$ is closed.}
\item {3) $0\cdot 0=0$.}
\item {4) $(\Bbb S,\cdot)$ has an identity element $1$ where $1\not=0$.}
\item {5) There exists $s\in\Bbb S$ such that $0\cdot s=1$.}
\item {6) $(\Bbb S,\cdot)$ is commutative.}
\item {7) $(\Bbb S,+,\cdot)$ is s-associative.}
\item {8) $(\Bbb S,+,\cdot)$ is wheel distributive.}
\end{list} Furthermore, we say that $(\Bbb S,+,\cdot)$ is an {\bf Ascended Field} if it satisfies the above axioms along with the additional axiom
\begin{center}
9) $(\Bbb S,+,\cdot)$ has scalar inverses.\,\,\,\,\,\,\,\,\,\,\,\,\,\,\,\,\,\,\,\,\,\,\,\,\,\,\,\,\,\,\,\,\,\,\,\,\,\,\,\,\,\,\,\,\,\,\,\,\,\,\,\,\,\,\,\,\,\,\,\,\,\,\,\,\,\,\,\,\,\,\,\,\,\,\,\,\,\,\,\,\,\,\,\,\,\,\,\,\,\,\,\,\,\,\,\,\,
\end{center}
By the first six axioms, we see that every ascended ring $({\Bbb S},+,\cdot)$ is also a commutative S-Structure.  This is true for every ascended field as well.  
  
\begin{proposition}
Let $(\Bbb S,+,\cdot)$ be an ascended ring (or an ascended field), then $0,1\in\Bbb S_0$ and $0,1\in \Lambda$ where $0$ and $1$ are distinct elements.
\end{proposition}

\begin{proof}
Since $1$ is the identity element of $(\Bbb S,\cdot)$, we know $1\in\Bbb S_0$ as $0\cdot 1=0$.  And since $0\cdot 0=0$, we know $0\in\Bbb S_0$.  Furthermore, since there exists $s\in {\Bbb S}$ such that $0\cdot s=1$, we know $s\in \Bbb S_1$.  This means that both $\Bbb S_0$ and $\Bbb S_1$ are non-empty, implying $0,1\in\Lambda$ as well.  Finally, $0$ and $1$ are distinct elements as $0\not=1$. $\square$
\end{proof}

\begin{proposition}
Let $(\Bbb S,+,\cdot)$ be an ascended ring (or an ascended field), then
\begin{list}{}{}
\item {1) If $x\in\Bbb S_0$ and $s,t\in\Bbb S$, then $x\cdot (s+t)=(x\cdot s)+(x\cdot t)$.}
\item {2) If $x\in\Bbb S_0$ and $s\in \Bbb S$, then $-(x\cdot s)=x\cdot (-s)$.}
\item {3) If $s_a\in\Bbb S_a$ and $s_b\in\Bbb S_b$, then $s_a+s_b=s_{a+b}$ where $s_{a+b}\in\Bbb S_{a+b}$.}
\end{list}
\end{proposition}

\begin{proof}
Let $x\in \Bbb S_0$ and let $s,t\in\Bbb S$.  This means that $x\cdot 0=0$, implying
$$x\cdot (s+t)=(x\cdot s)+(x\cdot t)-(x\cdot 0)=(x\cdot s)+(x\cdot t)-0=(x\cdot s)+(x\cdot t).$$
Therefore $x\cdot(s+t)=(x\cdot s)+(x\cdot t)$.  We can then see that
$$\begin{array}{rcl}
0&=&s + (-s) \cr
x\cdot 0&=&x\cdot (s+(-s)) \cr
0&=&(x\cdot s)+(x\cdot (-s)) \cr
-(x\cdot s)&=&x\cdot (-s). \cr
\end{array}$$
Therefore if $x\in \Bbb S_0$ and $s\in\Bbb S$, then $-(x\cdot s)=x\cdot (-s)$.
\newline \indent Finally, let $s_a\in\Bbb S_a$ and $s_b\in\Bbb S_b$, then we know that $0\cdot s_a=a$ and $0\cdot s_b=b$.  This means that $(0\cdot s_a)+(0\cdot s_b)=a+b$.  But since $0\in\Bbb S_0$, we see that
$$0\cdot (s_a+s_b)=(0\cdot s_a)+(0\cdot s_b)=a+b.$$ 
Therefore $0\cdot (s_a+s_b)=a+b$ implying $s_a+s_b=s_{a+b}$ where $s_{a+b}\in\Bbb S_{a+b}$. $\square$
\end{proof}

For each $\Bbb S_{\alpha}$ where $\alpha\in\Lambda$, we can fix an element $s_{\alpha}\in\Bbb S_{\alpha}$.  The following proposition will prove that for all $x\in\Bbb S_0$, elements of the form $s_{\alpha}+x$ also belong to $\Bbb S_{\alpha}$.  Furthermore, the proposition will prove that all other elements $t_{\alpha}\in\Bbb S_{\alpha}$ can be written in the form $s_{\alpha}+x_0$ for some $x_0\in\Bbb S_0$. 

\begin{proposition}
Let $(\Bbb S,+,\cdot)$ be an ascended ring (or an ascended field), then
\begin{list}{}{}
\item {1) If $s_a\in\Bbb S_a$, then $s_a+x\in\Bbb S_a$ for all $x\in\Bbb S_0$.}
\item {2) If $s_a,t_a\in\Bbb S_a$, then $t_a=s_a+x_0$ for some $x_0\in\Bbb S_0$.}
\end{list}
\end{proposition}

\begin{proof}
Let $s_a\in\Bbb S_a$ and let $x\in\Bbb S_0$, then $0\cdot x=0$.  Again, since $0\in\Bbb S_0$, this means that $0\cdot (s_a+x)=(0\cdot s_a)+(0\cdot x)=a+0=a$.  Therefore $0\cdot (s_a+x)=a$, which means that $s_a+x\in\Bbb S_a$ for all $x\in\Bbb S_0$.  
\newline \indent If we now let $t_a\in\Bbb S_a$, then $(0\cdot t_a)-(0\cdot s_a)=a-a=0$.  But at the same time $(0\cdot t_a)-(0\cdot s_a)=(0\cdot t_a)+(0\cdot (-s_a))=0\cdot(t_a-s_a)$, implying $0\cdot(t_a-s_a)=0$.  Therefore $t_a-s_a=x_0$ where $x_0\in\Bbb S_0$, which means that $t_a=s_a+x_0$ for some $x_0\in\Bbb S_0$. $\square$
\end{proof}
As a result of this proposition, we find that for any $\alpha\in\Lambda$, we can now describe every element in the set $\Bbb S_\alpha$ in terms of a fixed element $s_\alpha\in\Bbb S_\alpha$.

\begin{corollary}
Let $(\Bbb S,+,\cdot)$ be an ascended ring (or an ascended field) and let $\alpha\in\Lambda$.  If $s_\alpha\in\Bbb S_\alpha$ is a fixed element, then we can write
$$\Bbb S_\alpha=\{s_\alpha+x\,|\,x\in\Bbb S_0\}.$$
\end{corollary}

\begin{proof}
Let $\alpha\in\Lambda$ and fix $s_\alpha\in\Bbb S_\alpha$.  Now let $t\in\{s_\alpha+x\,|\,x\in\Bbb S_0\}$, then $t=s_\alpha+x\in\Bbb S_\alpha$.  At the same time, if $t\in\Bbb S_\alpha$, then $t=s_\alpha+x_0$ for some $x_0\in\Bbb S_0$.  This means that $t\in\{s_\alpha+x\,|\,x\in\Bbb S_0\}$.  Therefore $\Bbb S_\alpha=\{s_\alpha+x\,|\,x\in\Bbb S_0\}$. $\square$
\end{proof} 

\begin{proposition}
Let $(\Bbb S,+,\cdot)$ be an ascended ring (or an ascended field).  If $s,t\in\Bbb S$, then $s=t$ if and only if $s=s_x+a$ and $t=s_x+b$ for a fixed $s_x\in\Bbb S_x$ with $a,b\in\Bbb S_0$ where $a=b$.  
\end{proposition}

\begin{proof}
Let $s,t\in\Bbb S$ where $s=t$, then $0\cdot s=0\cdot t$ implying $s,t\in\Bbb S_x$ for some $x\in\Lambda$.  By fixing an element $s_x\in\Bbb S_x$, we see that $s=s_x+a$ and $t=s_x+b$ where $a,b\in\Bbb S_0$.  Since $s=t$, we know $s_x+a=s_x+b$ implying $a=b$.  The converse is also true.  That is, if $s=s_x+a$ while $t=s_x+b$ where $a=b$, then $s=s_x+a=s_x+b=t$.  Therefore, $s=t$ if and only if $s=s_x+a$ while $t=s_x+b$ for a fixed $s_x\in\Bbb S_x$ with $a,b\in\Bbb S_0$ where $a=b$. $\square$
\end{proof}

We can also show that every ascended ring/ascended field is complete regular and is an extension of a commutative ring or a field, respectively.

\begin{proposition}
Let $({\Bbb S},+,\cdot)$ be an ascended ring (or ascended field), then
\begin{list}{}{}
\item {1) If $m,n\in{\Bbb S}_0$ and $s\in{\Bbb S}$, then $m\cdot (n\cdot s)=n\cdot (m\cdot s)$.}
\item {2) If $x\in\Bbb S_0$ and $s_{\alpha}\in \Bbb S_{\alpha}$, then $x\cdot s_{\alpha}=s_{x\cdot \alpha}$ where $s_{x\cdot \alpha}\in\Bbb S_{x\cdot \alpha}$.}
\end{list}
\end{proposition}

\begin{proof}
Let $m,n\in{\Bbb S}_0$ and $s\in {\Bbb S}$, then we see that
$$\begin{array}{rcl}
m\cdot (n\cdot s)&=&(m\cdot n)\cdot s-\big(\big((m-1)\cdot (n-1)\big) \cdot (0\cdot s)\big) \cr
&=&(n\cdot m)\cdot s-\big(\big((n-1)\cdot (m-1)\big) \cdot (0\cdot s)\big) \cr
&=&n\cdot (m\cdot s). \cr
\end{array}$$
Therefore $m\cdot (n\cdot s)=n\cdot (m\cdot s)$.  Now let $x\in{\Bbb S}_0$ and $s_{\alpha}\in{\Bbb S}_{\alpha}$, then $0\cdot s_{\alpha}=\alpha$ implying $x\cdot (0\cdot s_{\alpha})=x\cdot\alpha$.  However, we know that $x\cdot(0\cdot s_{\alpha})=0\cdot(x\cdot s_{\alpha})$.  This means that $0\cdot(x\cdot s_{\alpha})=x\cdot\alpha$ and so $x\cdot s_{\alpha}=s_{x\cdot \alpha}$ where $s_{x\cdot \alpha}\in{\Bbb S}_{x\cdot \alpha}$. $\square$
\end{proof} 

\begin{proposition}
Let $({\Bbb S},+,\cdot)$ be an ascended ring (or an ascended field), then $\Lambda={\Bbb S}_0$ making $(\Bbb S,+,\cdot)$ complete regular.
\end{proposition}

\begin{proof}
Let $n\in\Bbb S_0$.  Since $1\in\Lambda$, we know $\exists s_1\in\Bbb S_1$.  Furthermore, we know that $n\cdot s_1\in \Bbb S_{n\cdot 1}$.  However $n\cdot 1=n$, implying that $n\cdot s_1\in\Bbb S_n$.  This means that $\Bbb S_n\not=\emptyset$ and so $n\in\Lambda$, implying $\Bbb S_0\subseteq \Lambda$.  
\newline \indent Now if $n\in\Lambda$, we know that $\exists s_n\in\Bbb S_n$ such that $0\cdot s_n=n$.  Since $1$ is the unity of $\Bbb S_0$, we know that $1\cdot 0=0$.  But since $0\in\Bbb S_0$, this means that $(-1)\cdot 0=-(1\cdot 0)=-0=0$, implying $(-1)\cdot 0=0$.  And so
$$\begin{array}{rcl}
n&=&0\cdot s_n \cr
&=&0\cdot (1\cdot s_n) \cr
&=&(0\cdot 1)\cdot s_n-\big(\big((0-1)\cdot (1-1)\big)\cdot(0\cdot s_n)\big) \cr
&=&0\cdot s_n-\big(\big((-1)\cdot 0\big)\cdot (0\cdot s_n)\big) \cr
&=&0\cdot s_n-\big(0\cdot (0\cdot s_n)\big) \cr
&=&n-(0\cdot n). \cr
\end{array}$$
Therefore $n=n-(0\cdot n)$, which means that $0\cdot n=0$.  And so $n\in\Bbb S_0$, implying $\Lambda\subseteq \Bbb S_0$.   This means that $\Lambda=\Bbb S_0$, proving that $(\Bbb S,+,\cdot)$ is complete regular. $\square$
\end{proof} 

\begin{proposition}
Let $(\Bbb S,+,\cdot)$ be an ascended ring (or an ascended field), then $(\Bbb S_0,+)$ is a commutative group.
\end{proposition}

\begin{proof}
Since $({\Bbb S},+)$ is a commutative group, we know that $({\Bbb S}_0,+)$ is associative and commutative.  Furthermore, we know that $0\in{\Bbb S}_0$.  If $x\in {\Bbb S}_0$, then $0\cdot x=0$ while 
$$0\cdot(-x)=-(0\cdot x)=-0=0.$$  
This proves that if $x\in\Bbb S_0$, then $-x\in\Bbb S_0$.  Now let $x,y\in\Bbb S_0$, then we know that 
$$0\cdot(x+y)=(0\cdot x)+(0\cdot y)=0+0=0.$$
This means that if $x,y\in\Bbb S_0$, then $x+y\in\Bbb S_0$ proving that $({\Bbb S}_0,+)$ is closed.  Therefore $(\Bbb S_0,+)$ forms a commutative group. $\square$
\end{proof}
If $(\Bbb S,+,\cdot)$ is an ascended ring, then we call $(\Bbb S_0,+,\cdot)$ the {\bf Scalar Ring of $(\Bbb S,+,\cdot)$}.  If $(\Bbb S,+,\cdot)$ is an ascended field, we call $(\Bbb S_0,+,\cdot)$ the {\bf Scalar Field of $(\Bbb S,+,\cdot)$}.  

\begin{theorem} 
$ $
\begin{list}{}{}
\item {1) If $(\Bbb S,+,\cdot)$ is an ascended ring, then $(\Bbb S_0,+,\cdot)$ is a commutative ring.}
\item {2) If $(\Bbb S,+,\cdot)$ is an ascended field, then $(\Bbb S_0,+,\cdot)$ is a field.}
\end{list}
\end{theorem}

\begin{proof}
For starters, we know that $(\Bbb S_0,+)$ is a commutative group.  Furthermore, since $(\Bbb S,+,\cdot)$ is commutative, we know that $(\Bbb S_0,\cdot)$ is commutative.  We also know that $1\in\Bbb S_0$ where $1$ is the unity of $(\Bbb S,+,\cdot)$, making $1$ the unity of $(\Bbb S_0,+,\cdot)$.  Since $x\cdot (s+t)=(x\cdot s)+(x\cdot t)$, if we let $s=y$ and $t=z$ where $x,y,z\in\Bbb S_0$, we see that $(\Bbb S_0,+,\cdot)$ is distributive.  Furthermore, since $x,y\in\Bbb S_0$, we know that $0\cdot (x\cdot y)=x\cdot(0\cdot y)=x\cdot 0=0$ implying $0\cdot (x\cdot y)=0$.  Therefore $x\cdot y\in\Bbb S_0$, proving that $(\Bbb S_0,\cdot)$ is closed.
\newline \indent Finally, if $x,y,z\in\Bbb S_0$, then $(x-1),(y-1)\in\Bbb S_0$ implying $(x-1)\cdot(y-1)\in\Bbb S_0$.  This means that $\big(\big((x-1)\cdot (y-1)\big)\cdot (0\cdot z)\big)=\big((x-1)\cdot(y-1)\big)\cdot 0=0$, hence
$$x\cdot (y\cdot z)=(x\cdot y)\cdot z-\big(\big((x-1)\cdot (y-1)\big)\cdot (0\cdot z)\big)=(x\cdot y)\cdot z-0=(x\cdot y)\cdot z.$$
Therefore $(\Bbb S_0,\cdot)$ is associative, proving that if $(\Bbb S,+,\cdot)$ is an ascended ring, then $(\Bbb S_0,+,\cdot)$ forms a commutative ring.  Furthermore, if $(\Bbb S,+,\cdot)$ has scalar inverses, then $(\Bbb S_0,+,\cdot)$ has multiplicative inverses for every $x\in\Bbb S_0$ where $x\not=0$.  So if $(\Bbb S,+,\cdot)$ is an ascended field, then we can conclude that $(\Bbb S_0,+,\cdot)$ is a field. $\square$
\end{proof}

Let's now explore examples of individual ascended fields.

\begin{example} 
\rm Let $(F,+,\cdot)$ be a field and let $\Bbb S=F\times F$ with elements $s,t\in\Bbb S$.  We then make the following definitions for addition and multiplication:
$$\begin{array}{rrcll}
s{\bf +}t&=&(x,y){\bf +}(u,v)&=&(x+u\,,\,y+v). \cr
s{\bf \cdot} t&=&(x,y){\bf \cdot}(u,v)&=&(xu+y+v-xv-yu,\,xv+yu). \cr
\end{array}$$  
\end{example} We can show that $(\Bbb S,\bf{+},{\bf \cdot})$ forms an ascended ring.  We can also show that $(\Bbb S_0,+,\cdot)=(F\times \{0\},+,\cdot)$ and that $(F\times \{0\},+,\cdot)$ is isomorphic to $(F,+_F,\cdot_F)$ with the isomorphism $\phi:F\to F\times \{0\}$ where $\phi(x)=(x,0)$.  Furthermore, we can show that each element $(x,0)\in\Bbb S_0$ where $x\not=0$ has an inverse element $(x^{-1},0)\in\Bbb S_0$.  Therefore $(\Bbb S,+,\cdot)$ is an ascended field.

\begin{example} 
\rm Let $(F,+,\cdot)$ be a field and let $\Bbb S=F\times F$ with elements $s,t\in\Bbb S$.  We then make the following definitions for addition and multiplication:
$$\begin{array}{rrcll}
s+t&=&(x,y)+(u,v)&=&(x+u\,,\,y+v). \cr
s\cdot t&=&(x,y)\cdot(u,v)&=&(xu+y+v-xv-yu,\,yv+xv+yu). \cr
\end{array}$$
\end{example} This structure also forms an ascended field with $(\Bbb S_0,+,\cdot)=(F\times \{0\},+,\cdot)$ where $(F\times \{0\},+,\cdot)$ is isomorphic to $(F,+,\cdot)$.  But this ascended field is different from the ascended field in Example 2.  This means that at least two different ascended fields containing the same $(\Bbb S_0,+,\cdot)$ exist.  And so ascended fields are not unique to the scalar field $(\Bbb S_0,+,\cdot)$ that they extend.  In Section 3 however, we will consider the quotient space of a particular equivalence relation on an ascended field which does uniquely extend the scalar field $(\Bbb S_0,+,\cdot)$.
  
\section{Tuples of Elements from the Dual of an Ascended Field}
Let $(\Bbb S,+,\cdot)$ be an ascended field and consider two copies of $(\Bbb S,+,\cdot)$.  We will refer to one copy as $(\Bbb S,+,\cdot)$, but refer to the other copy as $(\Bbb S^{-1},+,\cdot)$, called the {\bf Inverse Ascended Field of $(\Bbb S,+,\cdot)$}.  This gives us $(\Bbb S,+,\cdot)$ and $(\Bbb S^{-1},+,\cdot)$, both of which are the same algebraic structure.  So every $s\in\Bbb S$ has a unique counterpart $s^{-1}\in\Bbb S^{-1}$ and every counterpart $s^{-1}\in\Bbb S^{-1}$ can be traced back to a unique element $s\in\Bbb S$.  Furthermore, the addition and multiplication operations act the same on elements in $\Bbb S$ as they do on their counterparts in $\Bbb S^{-1}$.  This means that for every $s,t\in\Bbb S$ with respective counterparts $s^{-1},t^{-1}\in\Bbb S^{-1}$,
$$(s+t)^{-1}=s^{-1}+t^{-1}\,\,\,\,\,\,\,\,\,\,\,\,\,\,\,{\rm and}\,\,\,\,\,\,\,\,\,\,\,\,\,\,\,(s\cdot t)^{-1}=s^{-1}\cdot t^{-1}.$$
Hence the function $\phi:\Bbb S\rightarrow \Bbb S^{-1}$ where $\phi(s)=s^{-1}$ is an isomorphism with $\phi^{-1}(s^{-1})=s$.  So we can think of $^{-1}$ as an operation $^{-1}:\Bbb S\cup \Bbb S^{-1}\rightarrow \Bbb S\cup\Bbb S^{-1}$ such that $s^{-1}=\phi(s)$ for every $s\in\Bbb S$ and $(s^{-1})^{-1}=\phi^{-1}(s^{-1})=s$ for every $s^{-1}\in\Bbb S^{-1}$, making $^{-1}$ an involution.  We will refer to the elements $s$ and $s^{-1}$ as {\bf Abstract Inverses} of each other.
\newline \indent Since $(\Bbb S,+,\cdot)$ is an ascended field, if $s\in\Bbb S$, then we know $s\in\Bbb S_x$ for some $\Bbb S_x\subset \Bbb S$ where $x\in\Bbb S_0$.  Similarly, since $(\Bbb S^{-1},+,\cdot)$ is an exact copy of $(\Bbb S,+,\cdot)$, if $s^{-1}\in\Bbb S^{-1}$, then $s^{-1}\in\Bbb S_x^{-1}$ for some $\Bbb S_x^{-1}\subset\Bbb S^{-1}$ with $x\in\Bbb S_0$ where
$$\Bbb S_x^{-1}=\{s_x^{-1}\in\Bbb S^{-1}\,|\,0^{-1}\cdot s_x^{-1}=x^{-1}\}.$$
Note that if $s_x\in\Bbb S_x$ and $s_x^{-1}\in\Bbb S_x^{-1}$ for some $x\in\Bbb S_0$, then $0\cdot s_x=x$ with $x\in\Bbb S_0$ while $0^{-1}\cdot s_x^{-1}=x^{-1}$ with $x^{-1}\in\Bbb S_0^{-1}$.  
\newline \indent Finally, we can prove that if two ascended fields are equal, then their inverse ascended fields are also equal.

\begin{proposition}
Let $(\Bbb S_1,+_1,\cdot_1)$ and $(\Bbb S_2,+_2,\cdot_2)$ be two ascended fields such that $(\Bbb S_1,+_1,\cdot_1)=(\Bbb S_2,+_2,\cdot_2)$, then
$$(\Bbb S^{-1}_1,+_1,\cdot_1)=(\Bbb S^{-1}_2,+_2,\cdot_2).$$  
\end{proposition}

\begin{proof}
We know that the set $\Bbb S_1=\Bbb S_1^{-1}$ as they both contain the same elements.  Similarly, we know that $\Bbb S_2=\Bbb S_2^{-1}$ as they both contain the same elements.  Since we know that $(\Bbb S_1,+_1,\cdot_1)=(\Bbb S_2,+_2,\cdot_2)$, we know that $\Bbb S_1=\Bbb S_2$.  Therefore, we can conclude that $\Bbb S_1^{-1}=\Bbb S_2^{-1}$.  For convenience, we will set $\Bbb S=\Bbb S_1=\Bbb S_2$ and $\Bbb S^{-1}=\Bbb S_1^{-1}=\Bbb S_2^{-1}$. 
\newline \indent If $s,t\in\Bbb S$, since $(\Bbb S_1,+_1,\cdot_1)=(\Bbb S_2,+_2,\cdot_2)$, we know $s+_1 t=s+_2 t$.  So for every $s^{-1},t^{-1}\in\Bbb S^{-1}$, we see that $s^{-1}+_1 t^{-1}=(s+_1 t)^{-1}=(s+_2 t)^{-1}=s^{-1}+_2 t^{-1}$.  We also see that $s^{-1}\cdot_1 t^{-1}=(s\cdot_1 t)^{-1}=(s\cdot_2 t)^{-1}=s^{-1}\cdot_2 t^{-1}$.  Therefore $s^{-1}+_1 t^{-1}=s^{-1}+_2 t^{-1}$ and $s^{-1}\cdot_1 t^{-1}=s^{-1}\cdot_2 t^{-1}$ for every $s^{1},t^{-1}\in\Bbb S^{-1}$.  So we can conclude that $(\Bbb S_1^{-1},+_1,\cdot_1)=(\Bbb S_2^{-1},+_2,\cdot_2)$. $\square$
\end{proof}
\newpage
This proposition can also be generalized to ascended fields $(\Bbb S_1,+,\cdot)$ and $(\Bbb S_2,+,\cdot)$ which are isomorphic instead of equal.  If $(\Bbb S_1,+,\cdot)\cong (\Bbb S_2,+,\cdot)$, then there exists an isomorphism $f:\Bbb S_1\rightarrow \Bbb S_2$.  We can then define the function $\psi:\Bbb S_1^{-1}\rightarrow \Bbb S_2^{-1}$ where
$$\psi(s^{-1})=f(s)^{-1}.$$
Since $f$ is an isomorphism, we know $f$ is a bijection.  This means that $\psi$ must also be a bijection.  And if $s^{-1}, t^{-1}\in\Bbb S_1^{-1}$, then 
$$\psi(s^{-1}+t^{-1})=\psi((s+t)^{-1})=f(s+t)^{-1}=(f(s)+f(t))^{-1}=f(s)^{-1}+f(t)^{-1}$$
$$\psi(s^{-1}\cdot t^{-1})=\psi((s\cdot t)^{-1})=f(s\cdot t)^{-1}=(f(s)\cdot f(t))^{-1}=f(s)^{-1}\cdot f(t)^{-1}.$$
Therefore $\psi(s^{-1}+t^{-1})=f(s)^{-1}+f(t)^{-1}=\psi(s^{-1})+\psi(t^{-1})$ while at the same time $\psi(s^{-1}\cdot t^{-1})=f(s)^{-1}\cdot f(t)^{-1}=\psi(s^{-1})\cdot \psi(t^{-1})$.  So if two ascended fields are isomorphic, then their inverse ascended fields are also isomorphic.

\section*{The Dual of an Ascended Field} 
Let $(\Bbb S,+,\cdot)$ be an ascended field with the inverse ascended field $(\Bbb S^{-1},+,\cdot)$.  By assigning a {\bf Level} and a {\bf Value} to each element in $\Bbb S$ and $\Bbb S^{-1}$, we can differentiate between elements in the two structures.  If $s_x\in\Bbb S_x$ where $x\not=0$ is a non-scalar element in $\Bbb S$, then we define level$(s_x)=1$ and value$(s_x)=x$.  But if $s_x^{-1}\in\Bbb S_x^{-1}$ where $x\not=0$ is a non-scalar element in $\Bbb S^{-1}$, then we define level$(s_x^{-1})=-1$ and value$(s_x^{-1})=x$.  If $x\in\Bbb S_0$ with abstract inverse $x^{-1}\in\Bbb S_0^{-1}$, then we define level$(x)=$ level$(x^{-1})=0$ and value$(x)=$ value$(x^{-1})=x$.  As a consequence, every non-scalar element in $\Bbb S$ is no longer the same object as its corresponding abstract inverse in $\Bbb S^{-1}$ since they have different levels.  However, scalar elements in $\Bbb S$ and their corresponding abstract inverse in $\Bbb S^{-1}$ are still the same object as their levels are defined to be the same.  It's worth noting that this does not change how the operations of addition and multiplication act on elements in $\Bbb S$ and $\Bbb S^{-1}$.  The structures $(\Bbb S,+,\cdot)$ and $(\Bbb S^{-1},+,\cdot)$ remain isomorphic ascended fields.
\newline \indent We can then define $\Bbb S^*$, called the {\bf Dual of Ascended Field $(\Bbb S,+,\cdot)$}, as
$$\Bbb S^*=\Bbb S\cup\Bbb S^{-1}.$$
Since we no longer consider $\Bbb S$ and $\Bbb S^{-1}$ to be sets containing the same elements, we know $\Bbb S^*\not=\Bbb S$ and $\Bbb S^*\not=\Bbb S^{-1}$.  We also define the sets $\Bbb S_x^*\subset\Bbb S^*$ where $x\in\Bbb S_0$ as
$$\Bbb S_x^*=\Bbb S_x\cup\Bbb S_x^{-1}.$$
Since $s_x\in\Bbb S_x$ and its abstract inverse $s_x^{-1}\in\Bbb S_x^{-1}$ are no longer the same object when $x\not=0$, we know $\Bbb S_x^*\not=\Bbb S_x$ and $\Bbb S_x^*\not=\Bbb S_x^{-1}$ for every non-zero $x\in\Bbb S_0$.  But every scalar $x\in\Bbb S_0$ and its abstract inverse $x^{-1}\in\Bbb S_0^*$ are still the same object.  This means that $\Bbb S_0^*=\Bbb S_0^{-1}=\Bbb S_0$ and that every $x^{-1}\in\Bbb S_0^{-1}$ can be expressed as $x\in\Bbb S_0$.  We conclude that every distinct element in $\Bbb S^*$ is either a scalar $x\in\Bbb S_0$, a non-scalar $s\in\Bbb S$, or a non-scalar $s^{-1}\in\Bbb S^{-1}$.

\section*{Tuples of Elements from $\Bbb S^*$}
We now consider tuples with a finite amount of components consisting of non-scalar elements from the dual $\Bbb S^*$ of an ascended field $(\Bbb S,+,\cdot)$.  These tuples can be expressed as 
$$\vec{s}=(s_{x_1}^*,s_{x_2}^*,...\,,s_{x_n}^*)$$
where $s_{x_i}^*\in\Bbb S_{x_i}^*$ with non-zero $x_i\in\Bbb S_0$ for each $i\in\{1,2,...\,,n\}$.  We will define the {\bf Level} of $\vec{s}$ to be the sum of the levels of its components and will define the {\bf Value} of $\vec{s}$ to be the product of the values of its components.  This means that
$${\rm level}(\vec{s})=\sum_{i=1}^n{\rm level}(s_{x_i}^*)\,\,\,\,\,\,\,\,\,\,\,{\rm and}\,\,\,\,\,\,\,\,\,\,\,{\rm value}(\vec{s})=\prod_{i=1}^n{\rm value}(s_{x_i}^*)=\prod_{i=1}^n x_i$$
Note that every non-zero scalar $x$ can be the value of a tuple as each $x_i\in\Bbb S_0$ is non-zero.  Each integer number $n\in\Bbb Z$ can be the level of a tuple.  If level$(\vec{s})=n$ where $n\in\Bbb Z$, then we call $\vec{s}$ an {\bf n-Tuple}.  In this way, we can consider a scalar $x\in\Bbb S_0$ to be a $0$-tuple with value $x$, a non-scalar element $s_x\in\Bbb S_x$ to be a $1$-tuple with value $x$, and a non-scalar element $s_x^{-1}\in\Bbb S_x^{-1}$ to be a $(-1)$-tuple with value $x$.  For a general n-tuple, if $n>0$, then we call $\vec{s}$ an {\bf Ascended Tuple} and if $n<0$, then we call $\vec{s}$ a {\bf Descended Tuple}.  It's worth stating that there does not exist any tuple with a finite amount of components that has a value of $0$.  
\newline \indent If we consider two tuples $\vec{s}$ and $\vec{t}$, we say that 
\begin{center}
$\vec{s}=\vec{t}$ if and only if $\vec{s}$ and $\vec{t}$ have the same level and the same value.
\end{center}
We know that $\vec{s}$ has the same level and value as itself.  This means $\vec{s}=\vec{s}$.  If we let $\vec{s}=\vec{t}$, then $\vec{s}$ and $\vec{t}$ have the same level and value, implying that $\vec{t}=\vec{s}$.  Following the same logic, if $\vec{t}=\vec{s}$, then $\vec{s}=\vec{t}$.  So $\vec{s}=\vec{t}$ if and only if $\vec{t}=\vec{s}$.  Finally, if $\vec{r}$ is also a tuple where $\vec{s}=\vec{t}$ and $\vec{t}=\vec{r}$,  then the tuples $\vec{s}$, $\vec{t}$, and $\vec{r}$ all have the same level and value.  This means $\vec{s}=\vec{r}$.  Therefore, $=$ forms an equivalence relation.
\newline \indent As a consequence of this equivalence relation, we can express any $0$-tuple $\vec{s}$ with value $x\in\Bbb S_0\setminus\{0\}$ in {\bf Standard Form} by expressing it as
$$\vec{s}=x$$
since both $\vec{s}$ and scalar $x$ have level $0$ and value $x$.  Similarly, we can express any $1$-tuple $\vec{s}$ with value $x\in\Bbb S_0\setminus\{0\}$ in standard form by expressing it as 
$$\vec{s}=s_x$$ 
where $s_x\in\Bbb S_x$ as both $\vec{s}$ and $s_x$ have level $1$ and value $x$.  We can express any $(-1)$-tuple with value $x\in\Bbb S_0\setminus\{0\}$ in standard form by expressing it as
$$\vec{s}=s_x^{-1}$$
where $s_x^{-1}\in\Bbb S_x^{-1}$ as both $\vec{s}$ and $s_x^{-1}$ have level $-1$ and value $x$.  If $n\in\Bbb Z$ where $n>1$, then we can express any n-tuple with value $x\in\Bbb S_0\setminus\{0\}$ in standard form by expressing it as
$$\vec{s}=(s_x,s_1,t_1,...\,,r_1)$$ 
containing $n=|n|$ components where the first component $s_x\in\Bbb S_x$ while the remaining components $s_1, t_1,...\,, r_1\in{\Bbb S_1}$ as $\vec{s}$ and $(s_x,s_1,t_1,...\,,r_1)$ both have level $n$ (positive) and value $x$.  Finally, if $n<-1$, then we can express any n-tuple with value $x\in\Bbb S_0\setminus\{0\}$ in standard form by expressing it as
$$\vec{s}=(s_x^{-1},s_1^{-1},t_1^{-1},...\,,r_1^{-1})$$ 
containing $|n|$ components where the first component $s_x^{-1}\in\Bbb S_x^{-1}$ while the remaining components $s_1^{-1}, t_1^{-1}, ...\,,r_1^{-1}\in\Bbb S_1^{-1}$ as $\vec{s}$ and $(s_x^{-1},s_1^{-1},t_1^{-1},...\,,r_1^{-1})$ both have level $n$ (negative) and value $x$.  
\newline \indent While we can now express any tuple in a standard form, the standard form of a tuple is not necessarily unique.  For example, the tuples $(s_{x\cdot y},s_1)$ and $(s_{x\cdot y},s_1+1)$ where $s_{x\cdot y}\in\Bbb S_{x\cdot y}$ and $s_1\in\Bbb S_1$ are two different standard forms for the tuple $(s_x,s_y)$ where $s_x\in\Bbb S_x$ and $s_y\in\Bbb S_y$.  In light of this, we will define the sets $[n,x]$ such that
$$[n,x]=\{{\rm Standard\,\,forms\,\,for\,\,}\vec{s}\,\,|\,{\rm level\,of\,}\vec{s}=n,\,\,{\rm value\,of\,}\vec{s}=x\}.$$
We also define $[\Bbb Z,\Bbb S_0]=\{[n,x]\,|\,n\in\Bbb Z, x\in\Bbb S_0\setminus\{0\}\}$.  Note that $[0,x]=\{x\}$ as the scalar $x$ is the only standard form for a $0$-tuple with value $x$.  Also note that $[1,x]=\Bbb S_x$ as $\Bbb S_x$ contains all standard forms for a $1$-tuple with value $x$ while $[-1,x]=\Bbb S_x^{-1}$ as $\Bbb S_x^{-1}$ contains all standard forms for a $(-1)$-tuple with value $x$.  We also say that a set $[n,x]\in[\Bbb Z,\Bbb S_0]$ has a level of $n$ and a value of $x$.  So the singlet $\{x\}$ where $x\in\Bbb S_0\setminus\{0\}$ has level $0$ and value $x$, each $\Bbb S_x$ where $x\not=0$ has level $1$ and value $x$, and each $\Bbb S_x^{-1}$ where $x\not=0$ has level $-1$ and value $x$.
\newline \indent If $[n,x],[m,y]\in[\Bbb Z,\Bbb S_0]$ where $n=m$ and $x=y$, then clearly $\vec{s}\in[n,x]$ implies $\vec{s}\in[m,y]$ as ${\rm level}(\vec{s})=n=m$ and ${\rm value}(\vec{s})=x=y$.  A similar argument shows $\vec{s}\in[m,y]$ implies $\vec{s}\in[n,x]$.  Therefore $[n,x]=[m,y]$.  If instead $[n,x]=[m,y]$, then $\vec{s}\in[n,x]$ implies $\vec{s}\in[m,y]$.  That mean ${\rm level}(\vec{s})=n=m$ and ${\rm value}(\vec{s})=x=y$.  Therefore $m=n$ and $x=y$.  We conclude that if $[n,x],[m,y]\in[\Bbb Z,\Bbb S_0]$, then $[n,x]=[m,y]$ if and only if $n=m$ and $x=y$.  

\section*{Multiplying Tuples by a Scalar}
Let $(\Bbb S,+,\cdot)$ be an ascended field with the inverse ascended field $(\Bbb S^{-1},+,\cdot)$ and dual $\Bbb S^*$.  Recall that if $\alpha,\beta\in\Bbb S_0$ where $\alpha\not=0$, then we know $\alpha\cdot \beta\in\Bbb S_0$ is also a scalar as $(\Bbb S_0,+,\cdot)$ is a field.  Furthermore, if $s_x\in\Bbb S$ where $s_x\in\Bbb S_x$ with $x\not=0$, then $\alpha\cdot s_x=s_{\alpha\cdot x}\in\Bbb S_{\alpha\cdot x}$.  And if $s_x^{-1}\in\Bbb S$ where $s_x^{-1}\in\Bbb S_x^{-1}$ with $x\not=0$, then $\alpha^{-1}\cdot s_x^{-1}=s_{\alpha\cdot x}^{-1}\in\Bbb S_{\alpha\cdot x}^{-1}$ where $\alpha^{-1}\in\Bbb S_0^{-1}$ as $(\Bbb S,+,\cdot)$ and $(\Bbb S^{-1},+,\cdot)$ are isomorphic.  But $\alpha^{-1}$ and $\alpha$ are the same, so $\alpha\cdot s_x^{-1}=\alpha^{-1}\cdot s_x^{-1}=s_{\alpha\cdot x}^{-1}\in\Bbb S_{\alpha\cdot x}^{-1}$.  Putting this all together, we can formally express multiplication of an element $s_x^*\in\Bbb S_x^*$ where $x\not=0$ by a non-zero scalar $\alpha\in\Bbb S_0\setminus\{0\}$ as 
$$\alpha \cdot s_x^*=\left\{
        \begin{array}{cll}
						s_{\alpha\cdot x} & \quad s_x^*=s_x\in \Bbb S_x\setminus\{\Bbb S_0\} \\
						s_{\alpha\cdot x}^{-1} & \quad s_x^*=s_x^{-1}\in \Bbb S_x^{-1}\setminus\{\Bbb S_0\} \\
        \end{array} \right.$$
where $s_{\alpha\cdot x}\in\Bbb S_{\alpha\cdot x}$ and $s_{\alpha\cdot x}^{-1}\in\Bbb S_{\alpha\cdot x}^{-1}$.  So we see that $\alpha \cdot s_x^*=s_{\alpha\cdot x}^*\in\Bbb S_{\alpha \cdot x}^*$.  Therefore, we know how to multiply every element in $\Bbb S^*$ by a non-zero scalar. 
\newline \indent This scalar multiplication can be extended to any n-tuple $\vec{s}=(s_x^*, s_y^*, ...\,, s_z^*)$ containing two or more components where the components $s_x^*\in\Bbb S_x^*\setminus\{\Bbb S_0\}$, $s_y^*\in\Bbb S_y^*\setminus\{\Bbb S_0\}$, ...\,, $s_z^*\in\Bbb S_z^*\setminus\{\Bbb S_0\}$ by defining
$$\alpha\cdot (s_x^*, s_y^*, ...\,, s_z^*)=(\alpha \cdot s_x^*, s_y^*, ...\,, s_z^*)=(s_x^*, s_y^*, ...\,, s_z^*)\cdot\alpha.$$
Notice that only the first component is being scalar multiplied.  However, the next proposition shows that we don't have to choose only the first component.

\begin{proposition}
Let $\vec{s}=(s_x^*, s_y^*, ...\,, s_z^*)$ be an n-tuple containing two or more components where $s_x^*\in\Bbb S_x^*\setminus\{\Bbb S_0\}$, $s_y^*\in\Bbb S_y^*\setminus\{\Bbb S_0\}$, ...\,, $s_z^*\in\Bbb S_z^*\setminus\{\Bbb S_0\}$ and let $\alpha\in\Bbb S_0\setminus\{0\}$ be a scalar.  Then we have the following:
$$\begin{array}{rcl}
\alpha\cdot(s_x^*, s_y^*, ...\,, s_z^*)&=&(\alpha\cdot s_x^*, s_y^*, ...\,, s_z^*) \cr
&=&(s_x^*, \alpha\cdot s_y^*, ...\,, s_z^*) \cr
&& . \cr
&& . \cr
&=&(s_x^*, s_y^*, ...\,, \alpha\cdot s_z^*). \cr
\end{array}$$
\end{proposition}

\begin{proof} 
Let $\vec{s}=(s_x^*, s_y^*, ...\,, s_z^*)$ and let $\alpha\in\Bbb S_0\setminus\{0\}$.  We know that 
$$\begin{array}{rcl}
\alpha\cdot(s_x^*, s_y^*, ...\,, s_z^*)&=&(\alpha \cdot s_x^*, s_y^*, ...\,, s_z^*) \cr
&=&(s_{\alpha\cdot x}^*, s_y^*, ...\,, s_z^*). \cr
\end{array}$$
And because $(\Bbb S_0,+,\cdot)$ is a field with $\alpha, x,y,...\,,z\in\Bbb S_0$, we know that
$$\begin{array}{rcl}
\alpha \cdot x \cdot y \cdot ... \cdot z&=&(\alpha \cdot x) \cdot y \cdot ... \cdot z \cr
&=&x \cdot (\alpha \cdot y) \cdot ... \cdot z \cr
&& . \cr
&& . \cr
&=&x \cdot y \cdot ... \cdot (\alpha \cdot z). \cr
\end{array}$$
This means that the following n-tuples must be equal, as they all have the same level of $n$ and value of $\alpha\cdot x\cdot y\cdot ... \cdot z$:
$$\begin{array}{rcl}
\alpha\cdot(s_x^*, s_y^*, ...\,, s_z^*)&=&(s_{\alpha\cdot x}^*, s_y^*, ...\,, s_z^*) \cr
&=&(s_x^*, s_{\alpha\cdot y}^*, ...\,, s_z^*) \cr
&& . \cr
&& . \cr
&=&(s_x^*, s_y^*, ...\,, s_{\alpha\cdot z}^*). \cr
\end{array}$$
Finally, by replacing $s_{\alpha\cdot x}^*$ with $\alpha \cdot s_x^*$, replacing $s_{\alpha\cdot y}^*$ with $\alpha\cdot s_y^*$, ... , and replacing $s_{\alpha\cdot z}^*$ with $\alpha\cdot s_z^*$, we finish the proof. $\square$
\end{proof}
Multiplying a tuple by a scalar $\alpha$ doesn't change the level of the tuple, but does multiply the value of the tuple by $\alpha$.  So if $\vec{s}_x$ is an $n$-tuple with value $x\in\Bbb S_0\setminus\{0\}$ and $\alpha\in\Bbb S_0\setminus\{0\}$ is a scalar, then level$(\alpha\cdot \vec{s}_x)=n=0+n={\rm level}(\alpha)+{\rm level}(\vec{s}_x)$ and value$(\alpha\cdot \vec{s}_x)=\alpha\cdot{\rm value}(\vec{s}_x)={\rm value}(\alpha)\cdot{\rm value}(\vec{s}_x)$.

\section*{General Multiplication of Tuples}
Let $\vec{s}$ be an n-tuple with value $\alpha\in\Bbb S_0\setminus\{0\}$ and let $\vec{t}$ be an m-tuple with value $\beta\in\Bbb S_0\setminus\{0\}$.  If either $\vec{s}$ or $\vec{t}$ (or possibly both) are scalars, then we already know how to multiply $\vec{s}\cdot\vec{t}$.  If however both are not scalars, then we can use the Cartesian product to multiply $\vec{s}\cdot\vec{t}$.  In this case, we can express $\vec{s}$ and $\vec{t}$ as $\vec{s}=(s^*_x,s^*_y,...\,,s^*_z)$ and $\vec{t}=(s^*_u,s^*_v,...\,,s^*_w)$ where both tuples contain at least one component, $x\cdot y\cdot...\,\cdot z=\alpha$, and $u\cdot v\cdot...\,\cdot w=\beta$.  This means 
$$\begin{array}{rcl}
\vec{s}\cdot\vec{t}&=&(s^*_x,s^*_y,...\,,s^*_z)\cdot(s^*_u,s^*_v,...\,,s^*_w) \cr
&=&(s^*_x,s^*_y,...\,,s^*_z,s^*_u, s^*_v,...\,,s^*_w). \cr
\end{array}$$
We can then see that $\vec{s}\cdot\vec{t}$ has level $m+n$ and value 
$$x\cdot y\cdot...\cdot z\cdot u\cdot v\cdot...\cdot w=\alpha\cdot\beta.$$
Notice that level$(\vec{s}\cdot \vec{t})={\rm level}(\vec{s})+{\rm level}(\vec{t})$ and value$(\vec{s}\cdot \vec{t})={\rm value}(\vec{s})\cdot{\rm value}(\vec{t})$.  This is the same property we have when multiplying tuples with a scalar.  
\newline \indent Before moving on, note that if $\vec{s}$ is an $n$-tuple with value $x\in\Bbb S_0\setminus\{0\}$ and $\vec{t}$ is a $(-n)$-tuple with value $y\in\Bbb S_0\setminus\{0\}$, then $\vec{s}\cdot\vec{t}$ has level $n+(-n)=0$ and value $x\cdot y$.  So in standard form, $\vec{s}\cdot\vec{t}$ is just the non-zero scalar $x\cdot y$.	This applies to non-scalar elements $s_x\in\Bbb S_x$ and $s_y^{-1}\in\Bbb S_y^{-1}$ implying that $s_x\cdot s_y^{-1}=x\cdot y$.  Also note that multiplication of tuples, either with a scalar or in general, is a commutative operation. 

\section*{Addition of Tuples}
Consider $n$-tuple $\vec{s}_x$ and $m$-tuple $\vec{s}_y$ with values $x,y\in\Bbb S_0\setminus\{0\}$ respectively.  If $\vec{s}_x=x$ and $\vec{s}_y=y$, then $\vec{s}_x+\vec{s}_y=x+y$.  But this only results in a tuple if $x+y\not=0$.  If $\vec{s}_x$ and $\vec{s}_y$ are both non-scalars and $n\not=m$, then we define addition to result in the tuple with the higher level.  This means that
$$\vec{s}_x+\vec{s}_y=\left\{
        \begin{array}{cl}
						x+y & \quad \vec{s}_x=x,\,\vec{s}_y=y,\,x+y\not=0 \\
						\vec{s}_x & \quad {\rm level(}\vec{s}_x)> {\rm level(}\vec{s}_y) \\
						\vec{s}_y & \quad {\rm level(}\vec{s}_x)< {\rm level(}\vec{s}_y). \\
					\end{array} \right. $$	
If $s_x^*\in\Bbb S_x^*$ and $s_y^*\in\Bbb S_y^*$ are non-scalars where level$(s_x^*)=$ level$(s_y^*)$, then we see that $s_x^*+s_y^*=s_x+s_y=s_{x+y}\in\Bbb S_{x+y}$ or $s_x^*+s_y^*=s_x^{-1}+s_y^{-1}=s_{x+y}^{-1}\in\Bbb S_{x+y}^{-1}$.  In either case, $s_x^*+s_y^*=s_{x+y}^*\in\Bbb S_{x+y}^*$ where $s_x^*$, $s_y^*$, and $s_{x+y}^*$ have the same level.  So if $\vec{s}_x$ and $\vec{s}_y$ are both non-scalars with $n=m$ and $x+y\not=0$, then by expressing $\vec{s}_x$ and $\vec{s}_y$ in standard form as $\vec{s}_x=(s_x^*,s_1^*,...\,,t_1^*)$ and $\vec{s}_y=(s_y^*,s_1^*,...\,,t_1^*)$ where $s_x^*\in\Bbb S_x^*$ and $s_y^*\in\Bbb S_y^*$, using the same $s_1^*,...\,,t_1^*\in\Bbb S_1^*$, so that all of the components of $\vec{s}_x$ and $\vec{s}_y$ have the same level, we can define addition to be
$$\begin{array}{rcl}
\vec{s}_x+\vec{s}_y&=&(s_x^*,s_1^*,...\,,t_1^*)+(s_y^*,s_1^*,...\,,t_1^*) \cr
&=&(s_x^*+s_y^*,s_1^*,...\,,t_1^*). \cr
\end{array}$$  
In this way, $\vec{s}_x+\vec{s}_y$ is a tuple with level $n=m$ and value $x+y\not=0$ which is in standard form as $s_x^*+s_y^*=s_{x+y}^*\in\Bbb S_{x+y}^*$.  Note that this case of addition is not well-defined, but will still prove to be helpful.  Also note that because there are no tuples with value $0$, if $n=m$ and $x+y=0$, then $\vec{s}+\vec{t}$ does not exist.  

\section*{Infinite Tuples of Elements from $\Bbb S^*$}
Let $(\Bbb S,+,\cdot)$ be an ascended field with the inverse ascended field $(\Bbb S^{-1},+,\cdot)$ and dual $\Bbb S^*$.  We have been considering tuples containing non-scalar elements from $\Bbb S^*$ that consist of a finite amount of components.  We will now consider infinite tuples of elements from $\Bbb S^*$ which consist of a countably infinite amount of components.  Specifically, we will consider tuples of the form 
$$\vec{s}=(s_{x_1}^*,s_{x_2}^*,...\,,s_{x_n}^*,...)$$
where $s_{x_i}^*\in\Bbb S_{x_i}^*$ with non-zero $x_i\in\Bbb S_0$ for each $i\in\Bbb N$.  We know that
$${\rm level}(\vec{s})=\sum_{i=1}^\infty{\rm level}(s_{x_i}^*)\,\,\,\,\,\,\,\,\,\,\,{\rm and}\,\,\,\,\,\,\,\,\,\,\,{\rm value}(\vec{s})=\prod_{i=1}^\infty{\rm value}(s_{x_i}^*)=\prod_{i=1}^\infty x_i$$
Since level$(s_{x_i}^*)$ is either $1$ or $-1$, we will assume that $\sum_{i=1}^\infty{\rm level}(s_{x_i}^*)$ never converges.  If $\vec{s}$ satisfies $\sum_{i=1}^\infty{\rm level}(s_{x_i}^*)\rightarrow \infty$ in the traditional sense, that is the sequence of partial sums diverges to $\infty$, then we say that the level of $\vec{s}$ is $+\infty$.  If $\sum_{i=1}^\infty{\rm level}(s_{x_i}^*)\rightarrow -\infty$, again by considering the sequence of partial sums, then we say that the level of $\vec{s}$ is $-\infty$.  We will only consider infinite tuples whose levels are either $+\infty$ or $-\infty$.  At the same time, the product $\prod_{i=1}^\infty{\rm value}(s_{x_i}^*)$ stops making sense as we have not developed a way to evaluate it.  For this reason, we say that infinite tuples have no value.  
\newline \indent So if $\vec{s}$ and $\vec{t}$ are both infinite tuples, we define 
\begin{center}
$\vec{s}=\vec{t}$ if and only if they have the same level.
\end{center}
Clearly $\vec{s}=\vec{s}$ while $\vec{s}=\vec{t}$ if and only if $\vec{t}=\vec{s}$.  And if $\vec{r}$ is also an infinite tuple such that $\vec{s}=\vec{t}$ and $\vec{t}=\vec{r}$, then $\vec{s}=\vec{r}$ as they all have the same level.  So we can conclude that $=$ is an equivalence relation.  This allows us to express infinite tuples in standard form.  We can express any {\bf $(+\infty)$-tuple}, an infinite tuple with level $+\infty$, in standard form by expressing it as 
$$\vec{s}=(s_{x_1},s_{x_2},...\,,s_{x_n},...)$$
where $s_{x_i}\in\Bbb S_{x_i}$ with non-zero $x_i\in\Bbb S_0$ for each $i\in\Bbb N$.  We can express any {\bf $(-\infty)$-tuple}, an infinite tuple with level $-\infty$, in standard form by expressing it as 
$$\vec{s}=(s_{x_1}^{-1},s_{x_2}^{-1},...\,,s_{x_n}^{-1},...)$$
where $s_{x_i}^{-1}\in\Bbb S_{x_i}^{-1}$ with non-zero $x_i\in\Bbb S_0$ for each $i\in\Bbb N$.  The standard form of an infinite tuple is not unique and in light of this, we define the sets $0'$ and $\textcjheb{t}'$ such that
$$0'=\{{\rm Standard\,\,forms\,\,for\,\,}\vec{s}\,\,|\,{\rm level\,of\,}\vec{s}=-\infty\}\,$$
$$\textcjheb{t}'=\{{\rm Standard\,\,forms\,\,for\,\,}\vec{s}\,\,|\,{\rm level\,of\,}\vec{s}=+\infty\}.$$
We define the level of $0'$ to be $-\infty$ and the level of $\textcjheb{t}'$ to be $+\infty$.  

\section*{Addition with Infinite Tuples}
Let $\vec{s}$ be an $n$-tuple with $n\in\Bbb Z$, $\vec{s}_{+\infty}$ be a $(+\infty)$-tuple, and $\vec{s}_{-\infty}$ be a $(-\infty)$-tuple.  We see that 
$$\vec{s}_{-\infty}+\vec{s}=\vec{s}$$
$$\vec{s}_{-\infty}+\vec{s}_{+\infty}=\vec{s}_{+\infty}$$
as the level of any $n$-tuple and any $(+\infty)$-tuple will always be greater than the level of any $(-\infty)$-tuple.  We can also see that
$$\vec{s}_{+\infty}+\vec{s}=\vec{s}_{+\infty}$$
as the level of any $(+\infty)$-tuple is always greater than the level of any $n$-tuple.  
\newline \indent If $\vec{s}$ and $\vec{t}$ are two infinite tuples which have the same level, then they can be expressed in standard form as
$$\vec{s}=(s_x^*,s_u^*,...\,,s_v^*,...)\,\,\,\,\,\,\,{\rm and}\,\,\,\,\,\,\,\vec{t}=(s_y^*,s_u^*,...\,,s_v^*,...)$$ where $s_x^*\in\Bbb S_x^*$ and $s_y^*\in\Bbb S_y^*$, using the same $s_u^*,...\,,s_v^*,...\in\Bbb S_1^*$, we can define addition to be
$$\begin{array}{rcl}
\vec{s}+\vec{t}&=&(s_x^*,s_u^*,...\,,s_v^*,...)+(s_y^*,s_u^*,...\,,s_v^*,...) \cr
&=&(s_x^*+s_y^*,s_u^*,...\,,s_v^*,...). \cr
\end{array}$$  
In this way, $\vec{s}+\vec{t}$ is a tuple with the same level as both $\vec{s}$ and $\vec{t}$, either $+\infty$ or $-\infty$, and which is in standard form as $s_x^*+s_y^*=s_{x+y}^*\in\Bbb S_{x+y}^*$.  This mimics the addition of non-standard tuples with a finite amount of components with equal levels.  Also note that this case of addition is not well-defined, but will still prove to be helpful.  
\newline \indent  With that, we know how to add infinite tuples together, as well as add infinite tuples with n-tuples for any $n\in\Bbb Z$.  Note that these operations are commutative and that addition is well defined in the case of unequal levels while addition is not well defined in the case of equal levels.  However, in the case of equal levels, the sum of two $(+\infty)$-tuples is always a $(+\infty)$-tuple while the sum of two $(-\infty)$-tuples is always a $(-\infty)$-tuple.

\section*{Multiplication with Infinite Tuples}
We can express $n$-tuples with $n\in\Bbb Z$ containing one or more components as 
$$\vec{s}=(s_{x_1}^*,s_{x_2}^*,...\,,s_{x_n}^*)$$
where $s_{x_i}^*\in\Bbb S_{x_i}^*$ with non-zero $x_i\in\Bbb S_0$ for each $i\in\{1,2,...\,,n\}$.  We can express any $(+\infty)$-tuple as 
$$\vec{s}_{+\infty}=(s_{y_1}^*,s_{y_2}^*,...\,,s_{y_n}^*,...)$$
where $s_{y_i}^*\in\Bbb S_{y_i}^*$ with non-zero $y_i\in\Bbb S_0$ for each $i\in\Bbb N$ with $\sum_{i=1}^\infty{\rm level}(y_i)\rightarrow +\infty$.  Similarly, we can express any $(-\infty)$-tuple as 
$$\vec{s}_{-\infty}=(s_{z_1}^*,s_{z_2}^*,...\,,s_{z_n}^*,...)$$
where $s_{z_i}^*\in\Bbb S_{z_i}^*$ with non-zero $z_i\in\Bbb S_0$ for each $i\in\Bbb N$ with $\sum_{i=1}^\infty{\rm level}(z_i)\rightarrow -\infty$.  So if $\alpha\in\Bbb S_0$ is a non-zero scalar, then 
$$\begin{array}{rcl}
\alpha\cdot \vec{s}_{+\infty}&=& \alpha\cdot (s_{y_1}^*,s_{y_2}^*,...\,,s_{y_n}^*,...) \cr
&=&(\alpha\cdot s_{y_1}^*,s_{y_2}^*,...\,,s_{y_n}^*,...) \cr
\end{array}$$
while at the same time 
$$\begin{array}{rclcrcl}
\alpha\cdot \vec{s}_{-\infty}&=& \alpha\cdot (s_{z_1}^*,s_{z_2}^*,...\,,s_{z_n}^*,...) \cr
&=&(\alpha\cdot s_{z_1}^*,s_{z_2}^*,...\,,s_{z_n}^*,...).\cr
\end{array}$$
And since multiplying by a scalar does not change the level of a tuple, we know $\alpha\cdot\vec{s}_{+\infty}$ is a $(+\infty)$-tuple while $\alpha\cdot\vec{s}_{-\infty}$ is a $(-\infty)$-tuple.  We also see that
$$\begin{array}{rcl}
\vec{s}\cdot \vec{s}_{+\infty}&=&(s_{x_1}^*,s_{x_2}^*,...\,,s_{x_n}^*)\cdot (s_{y_1}^*,s_{y_2}^*,...\,,s_{y_n}^*,...) \cr
&=&(s_{x_1}^*,s_{x_2}^*,...\,,s_{x_n}^*,s_{y_1}^*,s_{y_2}^*,...\,,s_{y_n}^*,...) \cr
\end{array}$$
where ${\rm level}(\vec{s}\cdot\vec{s}_{+\infty})=\sum_{i=1}^n{\rm level}(s_{x_i})+\sum_{i=1}^{\infty}{\rm level}(s_{y_i})\rightarrow +\infty$ while 
$$\begin{array}{rcl}
\vec{s}\cdot \vec{s}_{-\infty}&=& (s_{x_1}^*,s_{x_2}^*,...\,,s_{x_n}^*)\cdot (s_{z_1}^*,s_{z_2}^*,...\,,s_{z_n}^*,...) \cr
&=&(s_{x_1}^*,s_{x_2}^*,...\,,s_{x_n}^*,s_{z_1}^*,s_{z_2}^*,...\,,s_{z_n}^*,...) \cr
\end{array}$$
where ${\rm level}(\vec{s}\cdot\vec{s}_{-\infty})=\sum_{i=1}^n{\rm level}(s_{x_i})+\sum_{i=1}^\infty{\rm level}(s_{z_i})\rightarrow -\infty$.  So we can conclude that $\vec{s}\cdot \vec{s}_{+\infty}$ is a $(+\infty)$-tuple while $\vec{s}\cdot \vec{s}_{-\infty}$ is a $(-\infty)$-tuple.  Finally, we see that
$$\begin{array}{rcl}
\vec{t}_{+\infty}\cdot \vec{s}_{+\infty}&=&(s_{x_1}^*,s_{x_2}^*,...\,,s_{x_n}^*,...)\cdot (s_{y_1}^*,s_{y_2}^*,...\,,s_{y_n}^*,...) \cr
&=&(s_{x_1}^*,s_{x_2}^*,...\,,s_{x_n}^*,...\,,s_{y_1}^*,s_{y_2}^*,...\,,s_{y_n}^*,...) \cr
\end{array}$$
where ${\rm level}(\vec{t}_{+\infty}\cdot \vec{s}_{+\infty})=\sum_{i=1}^{\infty}{\rm level}(s_{x_i})+\sum_{i=1}^{\infty}{\rm level}(s_{y_i})\rightarrow +\infty$ for every $(+\infty)$-tuple $s_{+\infty}$ and $t_{+\infty}$ while 
$$\begin{array}{rcl}
\vec{t}_{-\infty}\cdot \vec{s}_{-\infty}&=& (s_{w_1}^*,s_{w_2}^*,...\,,s_{w_n}^*,...)\cdot (s_{z_1}^*,s_{z_2}^*,...\,,s_{z_n}^*,...) \cr
&=&(s_{w_1}^*,s_{w_2}^*,...\,,s_{w_n}^*,...\,,s_{z_1}^*,s_{z_2}^*,...\,,s_{z_n}^*,...) \cr
\end{array}$$
where ${\rm level}(\vec{t}_{-\infty}\cdot \vec{s}_{-\infty})=\sum_{i=1}^{\infty}{\rm level}(s_{w_i})+\sum_{i=1}^\infty{\rm level}(s_{z_i})\rightarrow -\infty$ for every $(-\infty)$-tuple $s_{-\infty}$ and $t_{-\infty}$.  So we conclude that $\vec{t}_{+\infty}\cdot \vec{s}_{+\infty}$ is a $(+\infty)$-tuple while $\vec{t}_{-\infty}\cdot \vec{s}_{-\infty}$ is a $(-\infty)$-tuple.
\newline \indent However, in the case of $\vec{s}_{+\infty}\cdot\vec{s}_{-\infty}$, we see that
$$\begin{array}{rcl}
\vec{s}_{+\infty}\cdot \vec{s}_{-\infty}&=& (s_{y_1}^*,s_{y_2}^*,...\,,s_{y_n}^*,...)\cdot (s_{z_1}^*,s_{z_2}^*,...\,,s_{z_n}^*,...) \cr
&=&(s_{y_1}^*,s_{y_2}^*,...\,,s_{y_n}^*,...\,,s_{z_1}^*,s_{z_2}^*,...\,,s_{z_n}^*,...) \cr
\end{array}$$
where ${\rm level}(\vec{s}_{+\infty}\cdot\vec{s}_{-\infty})=\sum_{i=1}^\infty{\rm level}(y_i)+\sum_{i=1}^\infty{\rm level}(z_i)$ which may or may not exist.  Since we cannot guarantee the existence of the level of $\vec{s}_{+\infty}\cdot \vec{s}_{-\infty}$, we say that $\vec{s}_{+\infty}\cdot \vec{s}_{-\infty}$ is undefined.  It's worth noting that every other case of multiplication with infinite tuples is well defined and that multiplication with infinite tuples is a commutative operation, yielding the same results when the order in the operation is switched.


\newpage	
\section{Simple S-Extensions of Fields}
Let $({\Bbb S},+,\cdot)$ be an ascended field containing the scalar field $(\Bbb S_0,+,\cdot)$ with elements $x,y\in\Bbb S_0$.  In $(\Bbb S_0,+,\cdot)$, we know the relation $=$ satisfies
$$x=y\,\,{\rm if\,and\,only\,if}\,\,x-y=0.$$
Since $(\Bbb S_0,+,\cdot)$ is a field, this relation forms an equivalence relation where the equivalence class containing $x\in\Bbb S_0$ is just the singlet $\{x\}$ as
$$\begin{array}{rcl}
\{y\in\Bbb S_0\,|\,y=x\}&=&\{y\in\Bbb S_0\,|\,y-x=0\} \cr
&=&\{x\}. \cr
\end{array}$$
So for every $x\in\Bbb S_0$, there exists an equivalence class $\{x\}$.  This means that we can express the quotient set, which we denote as $\{\Bbb S_0\}$, as
$$\{\Bbb S_0\}=\{\{x\}\,|\,x\in\Bbb S_0\}.$$ 
Hence $=$ is an equivalence relation with the quotient set $\{\Bbb S_0\}=\{\{x\}\,|\,x\in\Bbb S_0\}$ containing the equivalence classes $\{x\}$.
\newline \indent By using the following addition and multiplication of sets:
$$A+B=\{a\,|\,a\in A\}+\{b\,|\,b\in B\}=\{a+b\,|\,a\in A, b\in B\}$$
$$A\cdot B=\{a\,|\,a\in A\}\cdot\{b\,|\,b\in B\}=\{a\cdot b\,|\,a\in A, b\in B\},$$
we can introduce addition and multiplication operations on $\{\Bbb S_0\}$ and find that for every $\{x\},\{y\}\in\{\Bbb S_0\}$, we have
$$\begin{array}{rcl}
\{x\}+\{y\}&=&\{x+y\} \cr
\{x\}\cdot \{y\}&=&\{x\cdot y\}. \cr
\end{array}$$
Under these addition and multiplication operations, $(\{\Bbb S_0\},+,\cdot)$ is isomorphic to the field $(\Bbb S_0,+,\cdot)$ as $\phi:\Bbb S_0\rightarrow \{\Bbb S_0\}$ where $\phi(x)=\{x\}$ is an isomorphism.  Hence we can think of the equivalence class $\{x\}\in\{\Bbb S_0\}$ as representing the scalar $x\in\Bbb S_0$.  In fact, for convenience sake, we will refer to $\{x\}$ as $x\in\Bbb S_0$ and only make the distinction between $\{x\}$ and $x$ when needed.  Similarly, we will refer to the set $\{\Bbb S_0\}$ as $\Bbb S_0$, only making the distinction when needed.
\newline \indent Now let $({\Bbb S},+,\cdot)$ be an ascended field where $s,t\in\Bbb S$ and define $\equiv'$ as follows:
$$s\equiv' t\,\,{\rm if\,and\,only\,if}\,\,0\cdot s=0\cdot t.$$
We can easily see that $s\equiv' s$ and that $s\equiv' t$ iff $t\equiv' s$.  Furthermore, if $s\equiv' t$ and $t\equiv' r$, then $0\cdot s=0\cdot t=0\cdot r$ implying $s\equiv' r$.  Therefore $\equiv'$ is an equivalence relation.  But because $\Lambda=\Bbb S_0$, for every $x\in\Bbb S_0$, there exists $t_x\in\Bbb S_x$ such that $0\cdot t_x=x$.  So we can express the equivalence class containing $t_x$ as $[x]$ where
$$[x]=\{s\in\Bbb S\,|\, s\equiv' t_x\}=\{s\in\Bbb S\,|\,0\cdot s=0\cdot t_x\}=\{s\in\Bbb S\,|\,0\cdot s=x\}=\Bbb S_x\not=\emptyset.$$
This means that for every $x\in\Bbb S_0$, the equivalence class $[x]=\Bbb S_x$.  And since $\{\Bbb S_\alpha\}_{\alpha\in\Lambda}$ forms a partition on $\Bbb S$ and $\Lambda=\Bbb S_0$, we know that $\{\Bbb S_x\}_{x\in\Bbb S_0}=\{[x]\}_{x\in\Bbb S_0}$ forms a partition on $\Bbb S$.  So the quotient set, which we will denote as $[\Bbb S_0']$, is 
$$[\Bbb S_0']=\{[x]\}_{x\in\Bbb S_0}=\{[x]\,|\,x\in\Bbb S_0\}.$$
Therefore $\equiv'$ is an equivalence relation with the quotient set $[\Bbb S_0']=\{[x]\,|\,x\in\Bbb S_0\}$ containing the equivalence classes $[x]=\Bbb S_x$.  
\newline \indent
By using the addition operation $A+B=\{a+b\,|\,a\in A, b\in B\}$, we also find that
$$\begin{array}{rcl}
[x]+[y]&=&\{s_x\,|\,s_x\in\Bbb S_x\}+\{s_y\,|\,s_y\in\Bbb S_y\} \cr
&=&\{s_x+s_y\,|\,s_x\in\Bbb S_x,\, s_y\in\Bbb S_y\}. \cr
\end{array}$$
But since $s_x\in\Bbb S_x$ and $s_y\in\Bbb S_y$, we know that $s_x+s_y=s_{x+y}\in\Bbb S_{x+y}$.  So if $s_x+s_y\in\{s_x+s_y\,|\,s_x\in\Bbb S_x,\, s_y\in\Bbb S_y\}$, then $s_x+s_y\in\Bbb S_{x+y}$.  Furthermore, if $s_{x+y}\in\Bbb S_{x+y}$, then because $(\Bbb S_0,+,\cdot)$ is a field, we know $x\in\Bbb S_0$ and $y\in\Bbb S_0$.  And since $\Lambda=\Bbb S_0$, we know $\exists s_x\in\Bbb S_x$.  We can then consider $s_y=s_{x+y}-s_x$, which we know is in $\Bbb S_{(x+y)-x}=\Bbb S_y$.  Then $s_{x+y}=s_x+s_y\in\{s_x+s_y\,|\,s_x\in\Bbb S_x,\, s_y\in\Bbb S_y\}$.  Putting this all together, we find that 
$$\begin{array}{rcl}
[x]+[y]&=&\{s_x+s_y\,|\,s_x\in\Bbb S_x,\, s_y\in\Bbb S_y\} \cr
&=&\Bbb S_{x+y} \cr
&=&[x+y]. \cr
\end{array}$$ 
Therefore $[x]+[y]=[x+y]$ for every $[x],[y]\in[\Bbb S_0']$.  We can then define subtraction as $[x]-[y]=[x]+[-y]=[x-y]$, which we see is well defined for every $[x],[y]\in[\Bbb S_0']$. 
\newline \indent Finally, let $(\Bbb S,+,\cdot)$ be an ascended field with $s,t\in \Bbb S$ and define $\equiv$ as follows:
$$s\equiv t\,\,{\rm if\,and\,only\,if}\,\,0\cdot s=0\cdot t\not=0\,\,{\rm or}\,\,s,t\in\Bbb S_0\,\,{\rm with}\,\,s=t.$$
It is easy to show that $\equiv$ forms an equivalence relation with equivalence classes $[x]=\{s\in\Bbb S\,|\,0\cdot s=x\}=\Bbb S_x$ for all $x\in\Bbb S_0\setminus\{0\}$ and $\{x\}=\{s\in\Bbb S_0\,|\,s=x\}$ for all $x\in\Bbb S_0$.  Notice that if $s_x\in\Bbb S_x$ where $x\not=0$, then $s_x\in [x]$ and if $s_x=x\in\Bbb S_0$, then $s_x\in\{x\}$.  For convenience, we will define
$$[\Bbb S_0]\,\,=\{\,[x]\,|\,x\in\Bbb S_0\setminus\{0\}\}.$$
Therefore $\equiv$ is an equivalence relation with the quotient set $\{\Bbb S_0\}\cup[\Bbb S_0]$ containing the equivalence classes $[x]=\Bbb S_x$ for all $x\in\Bbb S_0\setminus\{0\}$ and $\{x\}$ for all $x\in\Bbb S_0$.  
\newline \indent We know how the equivalence classes $\{x\},\{y\}\in\{\Bbb S_0\}$ add and multiply together.  We also know how equivalence classes $[x],[y]\in[\Bbb S_0]$ add together.  And by using the addition operation $A+B=\{a+b\,|\,a\in A, b\in B\}$, we find that if $\{x\}\in\{\Bbb S_0\}$ and $[y]\in\Bbb [\Bbb S_0]$, then 
$$\begin{array}{rcl}
\{x\}+[y]&=&\{x\}+\{s_y\,|\,s_y\in\Bbb S_y\} \cr
&=&\{x+s_y\,|\,s_y\in\Bbb S_y\}. \cr
$$\end{array}$$
But since $x\in\Bbb S_0$, we know that $x+s_y\in\Bbb S_y$.  So if $x+s_y\in\{x+s_y\,|\,s_y\in\Bbb S_y\}$, then $x+s_y\in\Bbb S_y$.  Furthermore, if $s_y'\in\Bbb S_y$, then we can consider $s_y=s_y'-x$ which is in $\Bbb S_y$ so that $s_y'=x+s_y\in\{x+s_y\,|\,s_y\in\Bbb S_y\}$.  Putting this all together, we see that
$$\begin{array}{rcl}
\{x\}+[y]&=&\{x+s_y\,|\,s_y\in\Bbb S_y\} \cr
&=&\Bbb S_y \cr
&=&[y]. \cr
\end{array}$$
Therefore $\{x\}+[y]=[y]$ for every $\{x\}\in\{\Bbb S_0\}$ and every $[y]\in[\Bbb S_0]$.  
\newline \indent Using the multiplication operation $A\cdot B=\{a\cdot b\,|\,a\in A, b\in B\}$, we also see that
$$\begin{array}{rcl}
\{x\}\cdot [y]&=&\{x\}\cdot \{s_y\,|\,s_y\in\Bbb S_y\} \cr
&=&\{x\cdot s_y\,|\,s_y\in\Bbb S_y\}. \cr
\end{array}$$
If $x\not=0$, then we know that $x\cdot s_y=s_{x\cdot y}\in\Bbb S_{x\cdot y}$.  So if $x\cdot s_y\in\{x\cdot s_y\,|\,s_y\in\Bbb S_y\}$, then $x\cdot s_y\in\Bbb S_{x\cdot y}$.  Furthermore if $s_{x\cdot y}\in\Bbb S_{x\cdot y}$, then because $(\Bbb S_0,+,\cdot)$ is a field, we know $x\in\Bbb S_0$ and $y\in\Bbb S_0$.  We can then consider $s_y=x^{-1}\cdot s_{x\cdot y}$ where $x^{-1}\in\Bbb S_0$ is the multiplicative inverse of $x$, which we know is in $\Bbb S_{x^{-1}\cdot(x\cdot y)}=\Bbb S_y$.  Then $s_{x\cdot y}=x\cdot s_y\in\{x\cdot s_y\,|\,s_y\in\Bbb S_y\}$.  Putting this all together, we find that
$$\begin{array}{rcl}
\{x\}\cdot [y]&=&\{x\cdot s_y\,|\,s_y\in\Bbb S_y\} \cr
&=&\Bbb S_{x\cdot y} \cr
&=&[x\cdot y]. \cr
\end{array}$$
Therefore $x\cdot[y]=[x\cdot y]$ for every $x\in\Bbb S_0$ where $x\not=0$ and every $[y]\in[\Bbb S_0]$.  Very importantly however, if $x=0$ and $[y]\in[\Bbb S_0]$, since $[y]=\Bbb S_y$, we see that 
$$\{0\}\cdot [y]=\{0\}\cdot\{s_y\,|\,s_y\in \Bbb S_y\}=\{0\cdot s_y\,|\, s_y\in\Bbb S_y\}=\{y\}.$$
This means that $\{0\}\cdot [y]=\{y\}$.  
\newline \indent With that, we have created the structure $(\{\Bbb S_0\}\cup[\Bbb S_0],+,\cdot)$.  Again however, since $(\{\Bbb S_0\},+,\cdot)$ is isomorphic to $(\Bbb S_0,+,\cdot)$, we will refer to singlets $\{x\}$ as the scalar $x$ and refer to the set $\{\Bbb S_0\}$ as the set $\Bbb S_0$.  Hence we get the structure $(\Bbb S_0\cup[\Bbb S_0],+,\cdot)$ called a {\bf Simple S-Extension of the Field $(\Bbb S_0,+,\cdot)$}, short for a {\bf Simple Ascended Extension of the Field $(\Bbb S_0,+,\cdot)$}.  So if $(\Bbb S_0\cup[\Bbb S_0],+,\cdot)$ is a simple s-extension of the scalar field $(\Bbb S_0,+,\cdot)$ of an ascended field $(\Bbb S,+,\cdot)$, and elements $s,t\in\Bbb S_0\cup[\Bbb S_0]$, then the simple s-extension has the addition operation 
$${s+t}=\left\{
        \begin{array}{cll}
            x+y & \quad s=x\in \Bbb S_0,\,&t=y\in\Bbb S_0 \\
						{\rm [}y{\rm ]} & \quad s=x\in \Bbb S_0,\,&t=[y]\in[\Bbb S_0] \\
						{\rm [}x+y{\rm ]}& \quad s=[x]\in [\Bbb S_0],\,&t=[y]\in [\Bbb S_0],\,\,\,\,\,x+y\not=0 \\
        \end{array} \right. $$
and the scalar multiplication operation 
$${s\cdot t}=\left\{
        \begin{array}{cll}
				    x\cdot y & \quad s=x\in \Bbb S_0,\,&t=y\in\Bbb S_0 \\
            {\rm [}x\cdot y{\rm ]} & \quad s=x\in \Bbb S_0\setminus\{0\},\,&t=[y]\in [\Bbb S_0] \\
						y & \quad s=0,\,&t=[y]\in [\Bbb S_0]. \\
						\end{array} \right. $$
It's worth noting that $0\cdot[y]=y$ for every $[y]\in[\Bbb S_0]$.
\newline \indent For every $x\in\Bbb S_0$, there exists $-x\in\Bbb S_0$ and for every $[x]\in[\Bbb S_0]$, there exists $-[x]\in[\Bbb S_0]$ where $-[x]=[-x]$.  This allows us to define subtraction in a simple s-extension $(\Bbb S_0\cup[\Bbb S_0],+,\cdot)$ as $s-t=s+(-t)$.  However, since $[0]\notin\Bbb S_0\cup[\Bbb S_0]$, subtraction $[x]-[y]$ is not defined if $x=y$.  In fact, if we consider subtraction as $s-t=r$ such that $s=t+r$, then subtraction of $[x]-[x]$ is indeterminate for every $[x]\in[\Bbb S_0]$ as $[x]-[x]=a$ satisfies $[x]=[x]+a$ for every $a\in\Bbb S_0$.
\newline \indent Now that we know how a simple s-extension works, we will see that unlike an ascended field, a simple s-extension of a field is unique to the field it extends. 

\begin{proposition}
Let $(\Bbb S_1,+_1,\cdot_1)$ and $(\Bbb S_2,+_2,\cdot_2)$ be ascended fields with the scalar fields $(\Bbb S_0^1,+_1,\cdot_1)$ and $(\Bbb S_0^2,+_2,\cdot_2)$, respectively.  If $(\Bbb S_0^1,+_1,\cdot_1)=(\Bbb S_0^2,+_2,\cdot_2)$, then $(\Bbb S_0^1\cup[\Bbb S_0^1],+_1,\cdot_1)=(\Bbb S_0^2\cup[\Bbb S_0^2],+_2,\cdot_2)$.  
\end{proposition}

\begin{proof}
Now $(\Bbb S_0^1,+_1,\cdot_1)=(\Bbb S_0^2,+_2,\cdot_2)$, which means that $\Bbb S_0^1=\Bbb S_0^2$.  This implies that $[\Bbb S_0^1]=\{[x]\,|\,x\in\Bbb S_0^1\}=\{[x]\,|\,x\in\Bbb S_0^2\}=[\Bbb S_0^2]$.  Therefore $\Bbb S_0^1=\Bbb S_0^2$ and $[\Bbb S_0^1]=[\Bbb S_0^2$], which means that $\Bbb S_0^1\cup[\Bbb S_0^1]=\Bbb S_0^2\cup[\Bbb S_0^2]$. 
\newline \indent If $x,y\in\Bbb S_0^1=\Bbb S_0^2$ and $[x],[y]\in[\Bbb S_0^1]=[\Bbb S_0^2]$, since $(\Bbb S_0^1,+_1,\cdot_1)=(\Bbb S_0^2,+_2,\cdot_2)$, $x+_1 y=x+_2 y$ and $[x]+_1 [y]=[x+_1 y]=[x+_2 y]=[x]+_2 [y]$.  We also see that $x+_1 [y]=[y]=x+_2 [y]$.  Therefore $s+_1 t=s+_2 t$ for every $s,t\in\Bbb S_0\cup[\Bbb S_0]$.
\newline \indent Furthermore, since $(\Bbb S_0^1,+_1,\cdot_1)=(\Bbb S_0^2,+_2,\cdot_2)$, we also see that $x\cdot_1 y=x\cdot_2 y$.  And if $x\not=0$, then $x\cdot_1 [y]=[x\cdot_1 y]=[x\cdot_2 y]=x\cdot_2 [y]$.  At the same time, $0\cdot_1 [y]=y=0\cdot_2 [y]$.  Therefore $s\cdot_1 t=s\cdot_2 t$ for every $s,t\in\Bbb S_0\cup[\Bbb S_0]$.  We can then conclude that $(\Bbb S_0^1\cup[\Bbb S_0^1],+_1,\cdot_1)=(\Bbb S_0^2\cup[\Bbb S_0^2],+_2,\cdot_2)$. $\square$
\end{proof}
So every ascended field containing the same scalar field $(\Bbb S_0,+,\cdot)$ gives rise to the same simple s-extension of $(\Bbb S_0,+,\cdot)$.  As a result, any field $(F,+_F,\cdot_F)$ can be extended into a unique simple s-extension of $(F\times\{0\},+,\cdot)$ where $(F\times\{0\},+,\cdot)$ is the scalar field of an ascended field and is isomorphic to $(F,+_F,\cdot_F)$.

\begin{corollary}
Let $(F,+_F,\cdot_F)$ be any field, then there exists a unique simple s-extension of $(F\times\{0\},+,\cdot)$ where $(F\times\{0\},+,\cdot)$ is the scalar field of an ascended field $(F\times F,+,\cdot)$ and is isomorphic to $(F,+_F,\cdot_F)$.
\end{corollary}

\begin{proof}
Since $(F,+_F,\cdot_F)$ is a field, we can use either Example 2 or Example 3 to turn this field into an ascended field $(F\times F,+,\cdot)$ where $(\Bbb S_0,+,\cdot)=(F\times\{0\},+,\cdot)$.  Now we know $(F\times\{0\},+,\cdot)$ is isomorphic to $(F,+_F,\cdot_F)$ as $\phi:F\to F\times \{0\}$ where $\phi(x)=(x,0)$ is an isomorphism.  From there, we use the equivalence relation $s\equiv t$ iff $0\cdot s=0\cdot t\not=0$ or $s,t\in\Bbb S_0$ with $s=t$ to create the quotient set $\{\Bbb S_0\}\cup[\Bbb S_0]=\{F\times\{0\}\}\cup[F\times\{0\}]$.  As discussed, we will refer to this quotient set as $\Bbb S_0\cup[\Bbb S_0]=(F\times\{0\})\cup[F\times\{0\}]$.  
\newline \indent This gives rise to the structure $(\Bbb S_0\cup[\Bbb S_0],+,\cdot)=((F\times\{0\})\cup[F\times\{0\}],+,\cdot)$, the simple s-extension of $(F\times\{0\},+,\cdot)$.  And since simple s-extensions are equal provided the ascended fields they come from contain the same $(\Bbb S_0,+,\cdot)$, we know that $((F\times\{0\})\cup[F\times\{0\}],+,\cdot)$ is the same no matter which ascended field with $(\Bbb S_0,+,\cdot)=(F\times\{0\},+,\cdot)$ we choose.  Therefore, the simple s-extension $((F\times\{0\})\cup[F\times\{0\}],+,\cdot)$ is the only simple s-extension of $(F\times\{0\},+,\cdot)$. $\square$
\end{proof}
In light of this corollary, $(F\cup[F],+,\cdot)$, the simple s-extension of the field $(F,+_F,\cdot_F)$, and $((F\times\{0\})\cup[F\times\{0\}],+,\cdot)$, the simple s-extension of the scalar field $(F\times\{0\},+,\cdot)$ of an ascended field $(F\times F,+,\cdot)$, are often used interchangeably.  
\newline \indent The last proposition can be generalized to two ascended fields $(\Bbb S_1,+_1,\cdot_1)$ and $(\Bbb S_2,+_2,\cdot_2)$ containing the scalar fields $(\Bbb S_0^1,+_1,\cdot_1)$ and $(\Bbb S_0^2,+_2,\cdot_2)$ which are isomorphic instead of equal.  Since $(\Bbb S_0^1,+_1,\cdot_1)\cong(\Bbb S_0^2,+_2,\cdot_2)$, there exists an isomorphism $f:\Bbb S_0^1\to\Bbb S_0^2$.  We can then define the function $\phi:\Bbb S_0^1\cup[\Bbb S_0^1]\to\Bbb S_0^2\cup[\Bbb S_0^2]$ such that 
$$\phi(s)=\left\{
        \begin{array}{ll}
            f(x) & s=x\in\Bbb S_0^1 \\
						{\rm [} f(x){\rm ]} & s=[x]\in [\Bbb S_0^1]. \\
        \end{array} \right. $$
Now the function $\phi$ is an isomorphism.  So if $(\Bbb S_1,+_1,\cdot_1)$ and $(\Bbb S_2,+_2,\cdot_2)$ are two ascended fields containing the scalar fields $(\Bbb S_0^1,+_1,\cdot_1)$ and $(\Bbb S_0^2,+_2,\cdot_2)$, respectively, where $(\Bbb S_0^1,+_1,\cdot_1)\cong(\Bbb S_0^2,+_2,\cdot_2)$, then $(\Bbb S_0^1\cup[\Bbb S_0^1],+_1,\cdot_1)\cong(\Bbb S_0^2\cup[\Bbb S_0^2],+_2,\cdot_2)$.

\section*{Division in a Simple S-Extension of a Field}
Let $(\Bbb S_0,+,\cdot)$ be a field and let $(\Bbb S_0\cup[\Bbb S_0],+,\cdot)$ be the simple s-extension of this field.  We define division in $(\Bbb S_0\cup[\Bbb S_0],+,\cdot)$ in the traditional manner: 
\newline If $s,t\in \Bbb S_0\cup[\Bbb S_0]$, then 
$${s\over t}=q\,\,{\rm such\,that}\,\, t\cdot q=s.$$
Using this definition, we can see what division in $(\Bbb S_0\cup[\Bbb S_0],+,\cdot)$ looks like. 

\begin{proposition}
Let $(\Bbb S_0\cup[\Bbb S_0],+,\cdot)$ be the simple s-extension of the field $(\Bbb S_0,+,\cdot)$.  If $s,t\in\Bbb S_0\cup[\Bbb S_0]$, then we see that
$${s\over t}=\left\{
        \begin{array}{cll}
            {x\cdot y^{-1}} & \quad s=x\in \Bbb S_0, & t=y\in\Bbb S_0\setminus\{0\} \\
						{x\cdot y^{-1}} & \quad s=[x]\in [\Bbb S_0], & t=[y]\in[\Bbb S_0] \\
						{\rm [}x\cdot y^{-1}\rm{]} & \quad s=[x]\in [\Bbb S_0], & t=y\in\Bbb S_0\setminus\{0\} \\
            0 & \quad s=y\in\Bbb S_0\setminus\{0\}, & t=[y]\in[\Bbb S_0]. \\
        \end{array} \right. $$
\end{proposition}

\begin{proof}
Let $x,y\in\Bbb S_0$ and $[x],[y]\in[\Bbb S_0]$.  If ${s\over t}={x\over y}=q$ where $y\not=0$, then $y\cdot q=x$.  Since $y\not=0$, we know that $y\cdot q\in\Bbb S_0$ only occurs when $q\in\Bbb S_0$.  Furthermore, since $x,y,q\in\Bbb S_0$ where $(\Bbb S_0,+,\cdot)$ is a field and $y\cdot q=x$, we know that $q=x\cdot y^{-1}$ is the only unique solution.  Therefore ${x\over y}=x\cdot y^{-1}$.  
\newline \indent If ${s\over t}={[x]\over [y]}=q$, then $[y]\cdot q=[x]$.  We know that $[y]\cdot q\in[\Bbb S_0]$ only occurs when $q\in\Bbb S_0$.  This means that $[x]=[y]\cdot q=[y\cdot q]$.  And so $x=y\cdot q$ where $x,y,q\in\Bbb S_0$ with $y\not=0$, implying again that $q=x\cdot y^{-1}$ is the only unique solution.  Therefore ${[x]\over [y]}=x\cdot y^{-1}$. 
\newline \indent If ${s\over t}={[x]\over y}=q$ where $y\not=0$, then $y\cdot q=[x]$.  Since $y\not=0$, we know that $y\cdot q\in[\Bbb S_0]$ only occurs when $q\in[\Bbb S_0]$.  So we can write $q$ as $q=[q_0]$ where $q_0\in\Bbb S_0$.  This means that $[x]=y\cdot q=y\cdot [q_0]=[y\cdot q_0]$, or $[x]=[y\cdot q_0]$.  And so $x=y\cdot q_0$ where $x,y,q_0\in\Bbb S_0$, implying that $q_0=y^{-1}\cdot x$ is the only unique solution.  Therefore ${[x]\over y}=[y^{-1}\cdot x]=[x\cdot y^{-1}]$.
\newline \indent Finally, if ${s\over t}={y\over [y]}=q$, then $q\cdot [y]=y$.  Furthermore, the only value of $q$ such that $q\cdot [y]=y$ is $q=0$.  Therefore ${y\over [y]}=0$. $\square$
\end{proof}
So we see that $\frac{s}{t}$ is defined if it takes the form $\frac{x}{y}$, $\frac{[x]}{[y]}$, $\frac{[x]}{y}$, or $\frac{x}{[x]}$.  However $\frac{x}{[x]}=0$ is the only division of the form ${x\over [y]}$ that is defined.  We can see this if we let $x\in\Bbb S_0$ and $[y]\in[\Bbb S_0]$ where $y\not=x$.  By letting ${x\over [y]}=q$, we see that $x=q\cdot [y]$.  This means that $q$ must be in $\Bbb S_0$.  But there is no element $q\in\Bbb S_0$ which satisfies this requirement as $q\cdot[y]=[q\cdot y]\notin\Bbb S_0$ if $q\not=0$ and $0\cdot [y]=y\not=x$ if $q=0$.  Therefore ${x\over [y]}$ does not exist if $y\not=x$. 
\newline \indent Before moving on, it is worth noting that the ratio of $x$ to $y$ is the same as the ratio of $[x]$ to $[y]$ as $\frac{[x]}{[y]}=\frac{x}{y}$.  It is also worth noting that if $s\in\Bbb S_0\cup[\Bbb S_0]$ and $y\in\Bbb S_0$ where $y\not=0$, then $\frac{s}{y}=y^{-1}\cdot s$ as $\frac{x}{y}=y^{-1}\cdot x$ and $\frac{[x]}{y}=[x\cdot y^{-1}]=y^{-1}\cdot[x]$.

\section*{Division By Zero in a Simple S-Extension of a Field}
Let $(\Bbb S_0,+,\cdot)$ be a field and let $(\Bbb S_0\cup[\Bbb S_0],+,\cdot)$ be the simple s-extension of this field.  If ${x\over 0}=q$ where $x\in\Bbb S_0$ is non-zero, then $0\cdot q=x$.  And we know that the only value of $q$ such that $0\cdot q=x$ with a non-zero $x$ is $q=[x]$.  Therefore, if $x\in\Bbb S_0$ where $x\not=0$, we can conclude that 
$${x\over 0}=[x].$$
This means that division by zero is well defined in the case of ${x\over 0}$ where $x\in\Bbb S_0$ is non-zero.  So if $(\Bbb S_0,+,\cdot)$ is a field and $(\Bbb S_0\cup[\Bbb S_0],+,\cdot)$ is the simple s-extension of this field with elements $x,y\in\Bbb S_0$, we can now write 
$${x\over y}= \left\{
        \begin{array}{cll}
            x\cdot y^{-1} & \quad y\not=0 \\
						{\rm [}x{\rm ]} & \quad y=0,\,x\not=0. \\
        \end{array} \right. $$
And since each $[x]$ is different for every choice of $x\not=0$, we see that every ${x\over 0}$ is an element unique to that particular $x$.  
\newline \indent We will finish up division by investigating ${0\over 0}$.  If ${0\over 0}=q$, then $0\cdot q=0$.  And we know that $0\cdot q=0$ is true for every $q\in\Bbb S_0$.  So we can conclude that ${0\over 0}\in\Bbb S_0$ is indeterminate.  Therefore, we have a unique algebraic structure in which division by zero satisfies a standard definition for the division operation and results in a different, unique element for ${x\over 0}$ for every non-zero $x$ in a field while satisfying our usual understanding of ${0\over 0}$.

\section*{Factorials of Negative Integers in the Simple S-Extension of $\Bbb R$}
Let $(\Bbb R\cup[\Bbb R],+,\cdot)$ be the simple s-extension of the real numbers $(\Bbb R,+,\cdot)$.  For $n\in\Bbb Z$, we define {\bf Factorials} by defining $0!=1$ and by using the recursive formula 
$$n!=n\cdot (n-1)!.$$
If we assume that $(n-1)!\in\Bbb R\cup[\Bbb R]$, we can calculate the $(n-1)^{th}$ factorial if we know the $n^{th}$ factorial.  

\begin{proposition} Let $(\Bbb R\cup[\Bbb R],+,\cdot)$ be the simple s-extension of the reals and let $n\in{\Bbb Z}$ where $n\not=0$.  If $(n-1)!\in\Bbb R\cup[\Bbb R]$, then
$$(n-1)!=\frac{n!}{n}.$$
\end{proposition}

\begin{proof} 
Let $n\in\Bbb Z$ be non-zero where $(n-1)!\in\Bbb R$ or $(n-1)!\in[\Bbb R]$.  If we let $(n-1)!=x\in\Bbb R$ or $(n-1)!=[x]\in[\Bbb R]$, then we see that $n\cdot(n-1)!=n\cdot x$ or $n\cdot(n-1)!=n\cdot [x]=[n\cdot x]$, respectively.  But by definition, we know that $n!=n\cdot (n-1)!$ which implies $n!=n\cdot x$ or $n!=[n\cdot x]$.  Furthermore, we see that either $n^{-1}\cdot n!=n^{-1}\cdot(n\cdot x)=(n^{-1}\cdot n)\cdot x=x=(n-1)!$ or $n^{-1}\cdot n!=n^{-1}\cdot [n\cdot x]=[n^{-1}\cdot (n\cdot x)]=[(n^{-1}\cdot n)\cdot x]=[x]=(n-1)!$. Either way, we see that $n^{-1}\cdot n!=(n-1)!$.  And since $n\in\Bbb Z$ is non-zero while $n!\in\Bbb R\cup[\Bbb R]$, we know that $n^{-1}\cdot n!=\frac{n!}{n}$.  So we can conclude that $(n-1)!=\frac{n!}{n}$. $\square$
\end{proof}
Now if $n=0$, then $n!=0!=1$.  This means that $n!=n\cdot(n-1)!$ becomes $1=0\cdot (-1)!$, implying that $(-1)!=[1]$ as no other element in $\Bbb R\cup[\Bbb R]$ satisfies this relationship.  Therefore $(-1)!=[1]$.  In fact, by using induction, we can prove that the following is a general formula for calculating $(-n)!$ for all negative integers $-n\in\Bbb Z$.

\begin{theorem}
Let $(\Bbb R\cup[\Bbb R],+,\cdot)$ be the simple s-extension of the reals and let $n\in\Bbb N$ where $n>0$, then
$$(-n)!=[{(-1)^{(n-1)}\over (n-1)!}].$$
\end{theorem}

\begin{proof}
If $n=1$, then $(-n)!=(-1)!=[1]$.  At the same time, we see that $[{(-1)^{(n-1)}\over (n-1)!}]=[{(-1)^0\over 0!}]=[\frac{1}{1}]=[1]$, implying $(-n)!=[{(-1)^{(n-1)}\over (n-1)!}]$ for $n=1$.  
\newline \indent If we assume $(-n)!=[{(-1)^{(n-1)}\over (n-1)!}]$ is true, then $[{(-1)^{(n-1)}\over (n-1)!}]=(-n)\cdot((-n)-1)!$.  This can only happen if $((-n)-1)!\in[\Bbb R]$.  Therefore
$$(-(n+1))!=((-n)-1)!=\frac{(-n)!}{-n}=\frac{[{(-1)^{(n-1)}\over (n-1)!}]}{-n}=[\frac{{(-1)^{(n-1)}\over (n-1)!}}{-n}]=[{(-1)^{n}\over n\cdot(n-1)!}]=[{(-1)^{n}\over n!}].$$
This means that if $(-n)!=[{(-1)^{(n-1)}\over (n-1)!}]$, then $(-(n+1))!=[{(-1)^{n}\over n!}]=[{(-1)^{((n+1)-1)}\over ((n+1)-1)!}]$.  So by induction, we can conclude that $(-n)!=[{(-1)^{(n-1)}\over (n-1)!}]$. $\square$
\end{proof}

\section*{Order in the Simple S-Extension of the Reals}
Let $(\Bbb R\cup[\Bbb R],+,\cdot)$ be the simple s-extension of the real numbers $(\Bbb R,+,\cdot)$.  Recall that the set of real numbers $\Bbb R$ contains a set of positive elements $\Bbb P(\Bbb R)$.  This means that the sets $[\Bbb R]$ and $[\Bbb R']$ also contain a set $[\Bbb P(\Bbb R)]$ such that
$$[\Bbb P(\Bbb R)]=\{[x]\,|\,x\in\Bbb P(\Bbb R)\}.$$  
Now, if $[x],[y]\in[\Bbb P(\Bbb R)]$, then $x,y\in\Bbb P(\Bbb R)$.  This means that $[x]+[y]=[x+y]$ where $x+y\in\Bbb P(\Bbb R)$.  Therefore $[x]+[y]\in[\Bbb P(\Bbb R)]$.  If $x\in\Bbb P(\Bbb R)$ and $[y]\in[\Bbb P(\Bbb R)]$, then $y\in\Bbb P(\Bbb R)$.  This means that $x\cdot [y]=[x\cdot y]$ where $x\cdot y\in\Bbb P(\Bbb R)$.  Therefore $x\cdot [y]\in[\Bbb P(\Bbb R)]$.  Furthermore, we see that if $[x]\in[\Bbb R']$, then either $x\in\Bbb P(\Bbb R)$, $-x\in\Bbb P(\Bbb R)$, or $x=0$.  This means that either $[x]\in[\Bbb P(\Bbb R)]$, $-[x]=[-x]\in[\Bbb P(\Bbb R)]$, or $[x]=[0]$.  We can then setup a natural ordering for elements $[x],[y]\in[\Bbb R']$ as
$$[x]<[y]\,\,{\rm if\,\,and\,\,only\,\,if}\,\,[y]-[x]\in[\Bbb P(\Bbb R)].$$
As a result, we see that
$$-[x]<[0]<[y]$$ 
for all positive $[x],[y]\in[\Bbb R']$.
\newline \indent We can now consider the set $\Bbb P(\Bbb R)\cup[\Bbb P(\Bbb R)]$.  If $x,y\in\Bbb P(\Bbb R)$ and $[x],[y]\in[\Bbb P(\Bbb R)]$, then $x+y\in\Bbb P(\Bbb R)$, $[x]+[y]\in[\Bbb P(\Bbb R)]$ and $[x]+y=[x]\in[\Bbb P(\Bbb R)]$.  Furthermore, we know that $x\cdot y\in\Bbb P(\Bbb R)$, implying that $x\cdot [y]=[x\cdot y]\in[\Bbb P(\Bbb R)]$.  Additionally, if $s\in\Bbb R\cup[\Bbb R]$, then either $s\in\Bbb P(\Bbb R)\cup[\Bbb P(\Bbb R)]$, $-s\in\Bbb P(\Bbb R)\cup[\Bbb P(\Bbb R)]$, or $s=0$.  Therefore, we can setup a natural ordering for elements $s,t\in\Bbb S_{\Bbb R}$ as
$$s<t\,\,{\rm if\,\,and\,\,only\,\,if}\,\,t-s\in\Bbb P(\Bbb R)\cup[\Bbb P(\Bbb R)].$$
Using this ordering, we find that for the simple s-extension $(\Bbb S_{\Bbb R},+,\cdot)$, we have
$$-[x]<-a<0<b<[y]$$
for all positive $[x],[y]\in[\Bbb R]$ and for all positive $a,b\in\Bbb R$.  With this in mind, we can consider any $[x]\in[\Bbb R]$ to be an {\bf Infinite Real Number}.  We can also consider any $x\in\Bbb S_0$ to be a {\bf Finite Real Number}. 

\section{Full and Complete S-Extensions of Fields}
Let $(\Bbb S,+,\cdot)$ be an ascended field with the inverse ascended field $(\Bbb S^{-1},+,\cdot)$.  To form a structure that extends simple s-extensions and has addition and multiplication operations defined for all elements in the structure, we will consider the non-traditional tuples from Section 2.  Recall from this section the set $[n,x]\in[\Bbb Z,\Bbb S_0]$ which is defined as
$$[n,x]=\{{\rm Standard\,\,forms\,\,for\,\,}\vec{s}\,\,\,|\,\,{\rm level\,of\,}\vec{s}=n,\,{\rm value\,of\,}\vec{s}=x\}$$
and is non-empty for every possible $n\in\Bbb Z$ and every possible $x\in\Bbb S_0\setminus\{0\}$.  
\newline \indent Also recall that $[0,x]=\{x\}$, $[1,x]=\Bbb S_x$, and $[-1,x]=\Bbb S_x^{-1}$ for any non-zero $x\in\Bbb S_0$.  At the same time, we know $[x]=\Bbb S_x$ in a simple s-extension.  This means that for every non-zero $x\in\Bbb S_0$, we have $[1,x]=\Bbb S_x=[x]$, or
$$[1,x]=[x].$$
Furthermore, we know that each singlet $\{x\}$ in a simple s-extension (which we usually refer to just as $x$) satisfies $[0,x]=\{x\}$ for every non-zero $x\in\Bbb S_0$.  Therefore, the set $[\Bbb Z,\Bbb S_0]$ of all sets $[n,x]$ contains every $[x]$ and every $\{x\}$ where $x\in\Bbb S_0$ is non-zero.  This means that 
$$(\{\Bbb S_0\}\cup[\Bbb S_0])\setminus\big\{\{0\}\big\}\subset[\Bbb Z,\Bbb S_0].$$
Hence $[\Bbb Z,\Bbb S_0]$ contains every non-zero element in the simple s-extension.  

\section*{Multiplication in $[\Bbb Z,\Bbb S_0]$}
We can introduce the following multiplication of sets:
$$A\cdot B=\{{\rm Standard\,\,forms\,\,for\,\,}a\cdot b\,\,|\,a\in A,\,\,b\in B\}.$$
If $[0,x],[0,y]\in[\Bbb Z,\Bbb S_0]$, then we see that
$$[0,x]\cdot[0,y]=\{x\}\cdot\{y\}=\{{\rm Standard\,\,forms\,\,for\,\,}x\cdot y\}=\{x\cdot y\}=[0,x\cdot y].$$
Therefore $[0,x]\cdot[0,y]=[0,x\cdot y]=[0+0,x\cdot y]$.  
\newline \indent At the same time, if $[m,y]\in[\Bbb Z,\Bbb S_0]$ where $m\not=0$, then
$$\begin{array}{rcl}
[0,x]\cdot[m,y]&=&\{x\}\cdot\{\,\vec{s}\,\,\,|\,\,\vec{s}\in[m,y]\} \cr
&=&\{{\rm Standard\,\,forms\,\,for\,\,}x\cdot \vec{s}\,\,\,|\,\,\vec{s}\in[m,y]\}. \cr
\end{array}$$
Since $\vec{s}\in[m,y]$, we know level$(\vec{s})=m$ and value$(\vec{s})=y$.  In Section 2, we saw that multiplying a tuple by a scalar $x$ does not change the level of the tuple, but multiplies the value of the tuple by $x$.  This means level$(x\cdot\vec{s})=m$ and value$(x\cdot\vec{s})=x\cdot y$.  So if $\vec{r}\in\{{\rm Standard\,\,forms\,\,for\,\,}x\cdot \vec{s}\,\,\,|\,\,\vec{s}\in[m,y]\}$, then $\vec{r}\in[m,x\cdot y]$.  Furthermore, if $\vec{t}\in[m,x\cdot y]$, then $\vec{t}$ can be expressed as
$$\begin{array}{rcl}
\vec{t}&=&(s_{x\cdot y}^*,s_1^*, ...\,, r_1^*) \cr
&=&(x\cdot s_y^*,s_1^*, ...\,, r_1^*) \cr
&=&x\cdot (s_y^*,s_1^*, ...\,, r_1^*). \cr
\end{array}$$
So we see that $\vec{t}=x\cdot (s_y^*,s_1^*, ...\,, r_1^*)$ where $(s_y^*,s_1^*, ...\,, r_1^*)\in[m,y]$.  Therefore $\vec{t}\in\{{\rm Standard\,\,forms\,\,for\,\,}x\cdot \vec{s}\,\,\,|\,\,\vec{s}\in[m,y]\}$.  Putting it all together, we find that
$$\begin{array}{rcl}
[0,x]\cdot[m,y]&=&\{{\rm Standard\,\,forms\,\,for\,\,}x\cdot \vec{s}\,\,\,|\,\,\vec{s}\in[m,y]\} \cr
&=&[m,x\cdot y]. \cr
\end{array}$$
So if $[0,x],[m,y]\in[\Bbb Z,\Bbb S_0]$ where $m\not=0$, then 
$$[0,x]\cdot [m,y]=[m,x\cdot y]=[0+m,x\cdot y].$$
Similarly, we can show that $[m,y]\cdot [0,x]=[m,y\cdot x]=[m+0,y\cdot x]$.  
\newline \indent Finally, if $[n,x],[m,y]\in[\Bbb Z,\Bbb S_0]$ where $n\not=0$ and $m\not=0$, then we see that
$$\begin{array}{rcl}
[n,x]\cdot[m,y]&=&\{\,\vec{s}\,\,\,|\,\,\vec{s}\in[n,x]\}\cdot\{\,\vec{t}\,\,\,|\,\,\vec{t}\in[m,y]\} \cr
&=&\{{\rm Standard\,\,forms\,\,for\,\,}\vec{s}\cdot \vec{t}\,\,\,|\,\,\vec{s}\in[n,x],\,\vec{t}\in[m,y]\}. \cr
\end{array}$$
In Section 2, we saw that if $\vec{s}$ is an n-tuple with value $x\in\Bbb S_0\setminus\{0\}$ and $\vec{t}$ is an m-tuple with value $y\in\Bbb S_0\setminus\{0\}$, then $s\cdot t$ is always a tuple with level $m+n$ and value $x\cdot y$.  So if $\vec{u}\in\{{\rm Standard\,\,forms\,\,for\,\,}\vec{s}\cdot \vec{t}\,\,\,|\,\,\vec{s}\in[n,x],\,\vec{t}\in[m,y]\}$, then $\vec{u}\in[n+m,x\cdot y]$.  Furthermore, if $\vec{r}\in[n+m,x\cdot y]$, then $\vec{r}$ is an $(n+m)$-tuple with value $x\cdot y\in\Bbb S_0\setminus\{0\}$.  So we can express $\vec{r}$ as 
$$\vec{r}=(s_x^*,s_y^*,s_1^*, t_1^*, ...\,, r_1^*,u_1^*)$$
where $s_x^*\in\Bbb S_x^*$, $s_y\in\Bbb S_y^*$, and $s_1^*, t_1^*, ...\,, r_1^*\in\Bbb S_1^*$.  Note that $(s_x^*,s_y^*,s_1^*, t_1^*, ...\,, r_1^*,u_1^*)$ is not necessarily in standard form.  So it has a level of $n+m$, but may not necessarily have $n+m$ components.  The components of $(s_x^*,s_y^*,s_1^*, t_1^*, ...\,, r_1^*,u_1^*)$ are chosen so that we can rewrite $\vec{r}$ as
$$\begin{array}{rcl}
\vec{r}&=&(s_x^*,s_y^*,s_1^*,t_1^*, ...\,, r_1^*,u_1^*) \cr
&=&(s_x^*,s_1^*, ...\,, r_1^*, s_y^*,t_1^*, ...\,,u_1^*) \cr
&=&(s_x^*,s_1^*, ...\,, r_1^*)\cdot(s_y^*,t_1^*, ...\,,u_1^*) \cr
\end{array}$$
where every component in $(s_x^*,s_1^*, ...\,, r_1^*)$ is either entirely from $\Bbb S$ or entirely from $\Bbb S^{-1}$, with the same case for the components of $(s_y^*,s_1^*, ...\,, r_1^*)$, while the level of $(s_x^*,s_1^*, ...\,, r_1^*)$ is $n$ and the level of $(s_y^*,t_1^*, ...\,, u_1^*)$ is $m$.  This means that 
$$\vec{r}=(s_x^*,s_1^*, ...\,, r_1^*)\cdot(s_y^*,t_1^*, ...\, ,u_1^*)$$
where $(s_x^*,s_1^*, ...\,, r_1^*)\in[n,x]$ and $(s_y^*,t_1^*, ...\,,u_1^*)\in[m,y]$.  Therefore, we see that $\vec{r}\in\{{\rm Standard\,\,forms\,\,for\,\,}\vec{s}\cdot \vec{t}\,\,\,|\,\,\vec{s}\in[n,x],\,\vec{t}\in[m,y]\}$.  Putting this all together, we see that  
$$\begin{array}{rcl}
[n,x]\cdot[m,y]&=&\{{\rm Standard\,\,forms\,\,for\,\,}\vec{s}\cdot \vec{t}\,\,\,|\,\,\vec{s}\in[n,x],\,\vec{t}\in[m,y]\} \cr
&=&[n+m,x\cdot y].
\end{array}$$
We conclude that if $[n,x],[m,y]\in[\Bbb Z,\Bbb S_0]$ with $n\not=0$ and $m\not=0$, then we know $[n,x]\cdot [m,y]=[n+m,x\cdot y]$.  
\newline \indent Hence we have shown that for every $[n,x],[m,y]\in[\Bbb Z,\Bbb S_0]$,
$$[n,x]\cdot[m,y]=[n+m,x\cdot y].$$
With this, we know how to multiply any two elements in $[\Bbb Z, \Bbb S_0]$ together.

\section*{Addition in $[\Bbb Z,\Bbb S_0]$}
We can also introduce the following addition of sets:
$$A+B=\{{\rm Standard\,\,forms\,\,for\,\,}a+b\,\,|\,a\in A,\,\,b\in B\}.$$
Let 
$[n,x],[m,y]\in[\Bbb Z,\Bbb S_0]$. 
If $n>m$, then we know that $\vec{s}+\vec{t}=\vec{s}$ for every $\vec{s}\in[n,x]$ and every $\vec{t}\in[m,y]$.  And since $\vec{s}\in[n,x]$, we know $\vec{s}$ is in standard form.  This means that
$$\begin{array}{rcl}
[n,x]+[m,y]&=&\{{\rm Standard\,\,forms\,\,for\,\,}\vec{s}+\vec{t}\,\,|\,\vec{s}\in [n,x],\,\,\vec{t}\in [m,y]\} \cr
&=&\{{\rm Standard\,\,forms\,\,for\,\,}\vec{s}\,\,|\,\vec{s}\in [n,x]\} \cr
&=&\{\,\vec{s}\,\,\,|\,\,\vec{s}\in [n,x]\} \cr
&=&[n,x]. \cr
\end{array}$$
So if $n>m$, then $[n,x]+[m,y]=[n,x]$. Similarly, if $n<m$, then we know $[n,x]+[m,y]=[m,y]$.  
\newline \indent If $n=m$ however, we know that 
$$\begin{array}{rcl}
[n,x]+[m,y]&=&[n,x]+[n,y] \cr
&=&\{{\rm Standard\,\,forms\,\,for\,\,}\vec{s}+\vec{t}\,\,|\,\vec{s}\in [n,x],\,\,\vec{t}\in [n,y]\}. \cr
\end{array}$$
In Section 2, we saw that $\vec{s}+\vec{t}$ is well defined if $n=0$ but is not well defined if $n\not=0$.  In either case though, if $x+y\not=0$, then $\vec{s}+\vec{t}$ has level $n$ and value $x+y$.  This means that if $\vec{u}\in\{{\rm Standard\,\,forms\,\,for\,\,}\vec{s}+\vec{t}\,\,|\,\vec{s}\in [n,x],\,\,\vec{t}\in [n,y]\}$, then $\vec{u}\in[n,x+y]$.  Furthermore, if $\vec{r}\in[n,x+y]$ where $x+y\not=0$, then we can write $\vec{r}$ as
$$\begin{array}{rcl}
\vec{r}&=&(s_{x+y}^*,s_1^*,...\,,t_1^*) \cr
&=&(s_{x}^*+s_{y}^*,s_1^*,...\,,t_1^*) \cr
&=&(s_{x}^*,s_1^*,...\,,t_1^*)+(s_{y}^*,s_1^*,...\,,t_1^*) \cr
\end{array}$$
where $s_{x+y}^*\in\Bbb S_{x+y}^*$ and $s_1^*,...\,,t_1^*\in\Bbb S_1^*$ with $s_{x+y}^*=s_x^*+s_y^*$ such that $s_x^*\in\Bbb S_x^*$ and $s_y^*\in\Bbb S_y^*$.  This means
$$\vec{r}=(s_{x}^*,s_1^*,...\,,t_1^*)+(s_{y}^*,s_1^*,...\,,t_1^*)$$
where $(s_{x}^*,s_1^*,...\,,t_1^*)\in[n,x]$ and $(s_{y}^*,s_1^*,...\,,t_1^*)\in[n,y]$.  Therefore, we see that
$\vec{r}\in\{{\rm Standard\,\,forms\,\,for\,\,}s+t\,\,|\,s\in [n,x],\,\,t\in [n,y]\}$.  Putting this all together, we see that if $x+y\not=0$, then
$$\begin{array}{rcl}
[n,x]+[m,y]&=&\{{\rm Standard\,\,forms\,\,for\,\,}s+t\,\,|\,s\in [n,x],\,\,t\in [n,y]\} \cr
&=&[n,x+y]. \cr
\end{array}$$
We conclude that if $[n,x],[m,y]\in[\Bbb Z,\Bbb S_0]$ with $n=m$ where $x+y\not=0$, then $[n,x]+[m,y]=[n,x+y]$.
\newline \indent Hence we have shown that for every $[n,x],[m,y]\in[\Bbb Z,\Bbb S_0]$,
$$[n,x]+[m,y]=\left\{
        \begin{array}{ll}
						{\rm [}n,x+y{\rm ]} & \quad {\rm if}\,\,n=m,\,\,x+y\not=0 \\
						{\rm [}n,x{\rm ]} & \quad {\rm if}\,\,n>m \\
						{\rm [}m,y{\rm ]} & \quad {\rm if}\,\,n<m. \\
        \end{array} \right.$$
With this, we know how to add elements in $[\Bbb Z,\Bbb S_0]$ together.

\section*{Full S-Extensions of Fields}
We call $([\Bbb Z,\Bbb S_0],+,\cdot)$ the {\bf Full S-Extension of the Field $(\Bbb S_0,+,\cdot)$}, with addition and multiplication operations
$$\begin{array}{rcl}
[n,x]+[m,y]&=&\left\{
				\begin{array}{cl}
						{\rm [}n,x+y{\rm ]} & \quad {\rm if}\,\,n=m,\,\,x+y\not=0 \\
						{\rm [}n,x{\rm ]} & \quad {\rm if}\,\,n>m \\
						{\rm [}m,y{\rm ]} & \quad {\rm if}\,\,n<m \\
        \end{array} \right. \cr
&& \cr
[n,x]\cdot [m,y]&=&[n+m, x\cdot y]. \cr				
\end{array}$$
Again, if $x\in\Bbb S_0\setminus\{0\}$, then $[0,x]=\{x\}$ and $[1,x]=[x]$ where $[x]\in[\Bbb S_0]$. We can show that the function $\phi:\{\Bbb S_0\}\cup[\Bbb S_0]\longrightarrow\{[n,x]\,\,|\,n\in\{0,1\},\,x\in\Bbb S_0\setminus\{0\}\}$ with $\phi(\{x\})=[0,x]$ and $\phi([x])=[1,x]$ is an isomorphism.  This means that a full s-extension extends the arithmetic of non-zero elements in a simple s-extension.  And just as elements $\{x\}\in\{\Bbb S_0\}$ in a simple s-extension are considered to be scalars in $\Bbb S_0$, we will consider elements $[0,x]\in[\Bbb Z,\Bbb S_0]$ to be scalars $x\in\Bbb S_0\setminus\{0\}$.  
\newline \indent We can then show that full s-extensions are unique to the field which they extend.  Furthermore, any field $(F,+_F,\cdot_F)$ can be extended into a unique full s-extension of the field $(F\times\{0\},+,\cdot)$, where $(F\times\{0\},+,\cdot)$ is the scalar field of an ascended field and is isomorphic to $(F,+_F,\cdot_F)$.  

\begin{proposition}
Let $(\Bbb S_1,+_1,\cdot_1)$ and $(\Bbb S_2,+_2,\cdot_2)$ be ascended fields with the scalar fields $(\Bbb S_0^1,+_1,\cdot_1)$ and $(\Bbb S_0^2,+_2,\cdot_2)$, respectively.  If $(\Bbb S_0^1,+_1,\cdot_1)=(\Bbb S_0^2,+_2,\cdot_2)$, then
$$([\Bbb Z,\Bbb S_0^1],+_1,\cdot_1)=([\Bbb Z, \Bbb S_0^2],+_2,\cdot_2).$$  
\end{proposition}

\begin{proof}
Now $(\Bbb S_0^1,+_1,\cdot_1)=(\Bbb S_0^2,+_2,\cdot_2)$, which means that $\Bbb S_0^1=\Bbb S_0^2$.  This implies $[\Bbb Z, \Bbb S_0^1]=\{[n,x]\,|\,n\in\Bbb Z, x\in\Bbb S_0^1\setminus\{0\}\}=\{[n,x]\,|\,n\in\Bbb Z, x\in\Bbb S_0^2\setminus\{0\}\}=[\Bbb Z, \Bbb S_0^2]$.  Therefore $[\Bbb Z, \Bbb S_0^1]=[\Bbb Z, \Bbb S_0^2]$. 
\newline \indent If $[n,x],[m,y]\in[\Bbb Z, \Bbb S_0^1]=[\Bbb Z, \Bbb S_0^2]=[\Bbb Z, \Bbb S_0]$, since $(\Bbb S_0^1,+_1,\cdot_1)=(\Bbb S_0^2,+_2,\cdot_2)$, we know that $x\cdot_1 y=x\cdot_2 y$ and that $x+_1 y=x+_2 y$.  This means that 
$$[n,x]\cdot_1 [m,y]=[n+m, x\cdot_1 y]=[n+m, x\cdot_2 y]=[n,x]\cdot_2 [m,y].$$
We also see that
$$\begin{array}{rcl}
[n,x]+_1 [m,y]&=&\left\{
				\begin{array}{cl}
						{\rm [}n,x+_1 y{\rm ]} & \quad {\rm if}\,\,n=m,\,\,x+y\not=0 \\
						{\rm [}n,x{\rm ]} & \quad {\rm if}\,\,n>m \\
						{\rm [}m,y{\rm ]} & \quad {\rm if}\,\,n<m \\
        \end{array} \right. \cr
				&&\cr
           &=&\left\{
				\begin{array}{cl}
						{\rm [}n,x+_2 y{\rm ]} & \quad {\rm if}\,\,n=m,\,\,x+y\not=0 \\
						{\rm [}n,x{\rm ]} & \quad {\rm if}\,\,n>m \\
						{\rm [}m,y{\rm ]} & \quad {\rm if}\,\,n<m \\
        \end{array} \right. \cr
				&&\cr
						&=&{\rm [}n,x{\rm ]}+_2 {\rm [}m,y{\rm ]}. \cr
\end{array}$$
Therefore $[n,x]\cdot_1 [m,y]=[n,x]\cdot_2 [m,y]$ and $[n,x]+_1 [m,y]=[n,x]+_2 [m,y]$ for every $[n,x],[m,y]\in[\Bbb Z, \Bbb S_0]$, implying $([\Bbb Z, \Bbb S_0^1],+_1,\cdot_1)=([\Bbb Z, \Bbb S_0^2],+_2,\cdot_2)$. $\square$
\end{proof}

\begin{corollary}
Let $(F,+_F,\cdot_F)$ be any field, then there exists a unique full s-extension of $(F\times\{0\},+,\cdot)$ where $(F\times\{0\},+,\cdot)$ is the scalar field of an ascended field $(F\times F,+,\cdot)$ and is isomorphic to $(F,+_F,\cdot_F)$.
\end{corollary}

\begin{proof}
Since $(F,+_F,\cdot_F)$ is a field, we can use either Example 2 or Example 3 to turn this field into an ascended field $(F\times F,+,\cdot)$ with $(\Bbb S_0,+,\cdot)=(F\times\{0\},+,\cdot)$ where $(F\times\{0\},+,\cdot)$ is isomorphic to $(F,+_F,\cdot_F)$.  This gives rise to the inverse ascended field $((F\times F)^{-1},+,\cdot)$ which we know is unique to $(F\times F,+,\cdot)$.  
\newline \indent From there, we can create the structure $([\Bbb Z, \Bbb S_0],+,\cdot)=([\Bbb Z, F\times\{0\}],+,\cdot)$, the full s-extension of $(F\times\{0\},+,\cdot)$, by considering the non-traditional tuples whose components are elements from $(F\times F)^*=(F\times F)\cup(F\times F)^{-1}$.  Since full s-extensions are equal provided the ascended fields they come from contain the same $(\Bbb S_0,+,\cdot)$, we know that $((\Bbb Z, F\times\{0\}),+,\cdot)$ is the same no matter which ascended field containing the field $(\Bbb S_0,+,\cdot)=(F\times\{0\},+,\cdot)$ we choose.  Therefore $([\Bbb Z, F\times\{0\}],+,\cdot)$ is the unique full s-extension of $(F\times\{0\},+,\cdot)$. $\square$
\end{proof}
In light of this corollary, $([\Bbb Z,F],+,\cdot)$, the full s-extension of the field $(F,+_F,\cdot_F)$, and $([\Bbb Z,F\times\{0\}],+,\cdot)$, the full s-extension of the scalar field $(F\times\{0\},+,\cdot)$ of an ascended field $(F\times F,+,\cdot)$, are often used interchangeably.  
\newline \indent This proposition can also be generalized to two ascended fields $(\Bbb S_1,+_1,\cdot_1)$ and $(\Bbb S_2,+_2,\cdot_2)$ containing the fields $(\Bbb S_0^1,+_1,\cdot_1)$ and $(\Bbb S_0^2,+_2,\cdot_2)$ which are isomorphic instead of equal.  Since $(\Bbb S_0^1,+_1,\cdot_1)\cong(\Bbb S_0^2,+_2,\cdot_2)$, there exists an isomorphism $f:\Bbb S_0^1\to\Bbb S_0^2$.  We can then define the function $\phi:(\Bbb Z, \Bbb S_0^1)\to(\Bbb Z, \Bbb S_0^2)$ such that 
\begin{center}
$\phi(s)=[n,f(x)]$ for all $s=[n,x]\in[\Bbb Z, \Bbb S_0^1]$.
\end{center}  
Now the function $\phi$ is an isomorphism.  So if $(\Bbb S_1,+_1,\cdot_1)$ and $(\Bbb S_2,+_2,\cdot_2)$ are two ascended fields containing the scalar fields $(\Bbb S_0^1,+_1,\cdot_1)$ and $(\Bbb S_0^2,+_2,\cdot_2)$, respectively, where $(\Bbb S_1,+_1,\cdot_1)\cong(\Bbb S_2,+_2,\cdot_2)$, then $([\Bbb Z, \Bbb S_0^1],+_1,\cdot_1)\cong([\Bbb Z, \Bbb S_0^2],+_2,\cdot_2)$. 
\newline \indent Finally, in a full s-extension $[\Bbb Z,\Bbb S_0]$, we know that for every $[n,x]\in[\Bbb Z,\Bbb S_0]$, there exists $-[n,x]=[n,-x]\in[\Bbb Z,\Bbb S_0]$.  So we can define a natural subtraction as $[n,x]-[m,y]=[n,x]+[m,-y]$.  Furthermore, this subtraction is similar to subtraction in simple s-extensions as $[n,x]-[n,x]$ is indeterminate for every $[n,x]\in[\Bbb Z,\Bbb S_0]$.  This is because any element $[m,y]\in[\Bbb Z,\Bbb S_0]$ where $m<n$ satisfies $[n,x]+[m,y]=[n,x]$.
 
\section*{Division in a Full S-Extension of a Field}
Let $(\Bbb S_0,+,\cdot)$ be a field and let $([\Bbb Z, \Bbb S_0],+,\cdot)$ be the full s-extension of the field.  Let $[n,x],[m,y]\in[\Bbb Z,\Bbb S_0]$, then we can define 
$$\frac{[n,x]}{[m,y]}=q\,\,\,{\rm such\,\,that}\,\,\,[m,y]\cdot q=[n,x].$$
If we let $q=[p,z]\in[\Bbb Z,\Bbb S_0]$ such that $\frac{[n,x]}{[m,y]}=[p,z]$, then $[m,y]\cdot [p,z]=[n,x]$.  This means that
$$[m+p,\,y\cdot z]=[n,x].$$
Therefore $n=m+p$ and $x=y\cdot z$, implying that $p=n-m$ and that $z=\frac{x}{y}$.  So using this definition, we see that
$$\frac{[n,x]}{[m,y]}=[n-m,\frac{x}{y}].$$
In simple s-extensions, division is defined for only a few cases and is not defined for all elements.  However in full s-extensions, division is defined for all elements $[n,x],[m,y]\in[\Bbb Z,\Bbb S_0]$.  Furthermore, $\frac{[n,x]}{[n,y]}=[n-n,\frac{x}{y}]=[0,\frac{x}{y}]$ which means the ratio of $[n,x]$ to $[n,y]$ is the same as the ratio of $x$ to $y$ for every $n\in\Bbb Z$.

\section*{Addition with Absolute Zero}
Let $(\Bbb S_0,+,\cdot)$ be a field and let $([\Bbb Z,\Bbb S_0],+,\cdot)$ be the full s-extension of this field.  Also recall the set $0'=\{{\rm Standard\,\,forms\,\,for\,\,}\vec{s}\,\,|\,{\rm level\,of\,}\vec{s}=-\infty\}$ from Section 2.  Now consider the set $[\Bbb Z,\Bbb S_0]^*$ where 
$$[\Bbb Z,\Bbb S_0]^*=\,[\Bbb Z,\Bbb S_0]\cup\{0'\}.$$
We then see that
$$0'+\,0'=\{{\rm Standard\,\,forms\,\,for\,\,}\vec{s}+\vec{t}\,\,|\,\vec{s}\in 0',\,\,\vec{t}\in 0'\}.$$
In Section 2, we saw that $\vec{s}+\vec{t}$ is also a $(-\infty)$-tuple.  This means that if $\vec{u}\in\{{\rm Standard\,\,forms\,\,for\,\,}\vec{s}+\vec{t}\,\,|\,\vec{s}\in 0',\,\,\vec{t}\in 0'\}$, then $\vec{u}\in 0'$.  Furthermore, if $\vec{r}\in 0'$, then we can write $\vec{r}$ as 
$$\begin{array}{rcl}
\vec{r}&=&(s_x^*,s_y^*,...\,,s_z^*,...) \cr
&=&(s_{x-1+1}^*,s_y^*,...\,,s_z^*,...) \cr
&=&(s_{x-1}^*+s_1^*,s_y^*,...\,,s_z^*,...) \cr
&=&(s_{x-1}^*,s_y^*,...\,,s_z^*,...)+(s_1^*,s_y^*,...\,,s_z^*,...). \cr
\end{array}$$
This means that 
$$\vec{r}=(s_{x-1}^*,s_y^*,...\,,s_z^*,...)+(s_1^*,s_y^*,...\,,s_z^*,...)$$
where $(s_{x-1}^*,s_y^*,...\,,s_z^*,...)\in 0'$ and $(s_1^*,s_y^*,...\,,s_z^*,...)\in 0'$.  Therefore, we see that $\vec{r}\in\{{\rm Standard\,\,forms\,\,for\,\,}\vec{s}+\vec{t}\,\,|\,\vec{s}\in 0',\,\,\vec{t}\in 0'\}$.  Putting this all together, we see that 
$$\begin{array}{rcl}
0'+\,0'&=&\{{\rm Standard\,\,forms\,\,for\,\,}\vec{s}+\vec{t}\,\,|\,\vec{s}\in 0',\,\,\vec{t}\in 0'\} \cr
&=&0'. \cr
\end{array}$$
We conclude that $0'+\,0'=0'$.
\newline \indent If however $[n,x]\in[\Bbb Z,\Bbb S_0]$, then from Section 2, we know that $\vec{s}+\vec{t}=\vec{s}$ for every $\vec{s}\in[n,x]$ and every $\vec{t}\in 0'$.  So we see that
$$\begin{array}{rcl}
[n,x]+\,0'&=&\{{\rm Standard\,\,forms\,\,for\,\,}\vec{s}+\vec{t}\,\,\,|\,\,\vec{s}\in[n,x],\,\,\vec{t}\in 0'\} \cr
&=&\{{\rm Standard\,\,forms\,\,for\,\,}\vec{s}\,\,\,|\,\,\vec{s}\in[n,x]\} \cr
&=&\{\,\vec{s}\,\,\,|\,\,\vec{s}\in[n,x]\} \cr
&=&[n,x]. \cr
\end{array}$$
Therefore $[n,x]+0'=[n,x]$.  Similarly, we can show that $0'+[n,x]=[n,x]$.  We conclude that $[n,x]+0'=0'+[n,x]=[n,x]$ for every $[n,x]\in[\Bbb Z,\Bbb S_0]$.  We also conclude that $0'$ acts as an additive identity for elements in $[\Bbb Z,\Bbb S_0]^*$. 

\section*{Multiplication with Absolute Zero}
We can also see that
$$0'\cdot \,0'=\{{\rm Standard\,\,forms\,\,for\,\,}\vec{s}\cdot\vec{t}\,\,|\,\vec{s}\in 0',\,\,\vec{t}\in 0'\}.$$
Now in Section 2, we saw that $\vec{s}\cdot \vec{t}$ is also a $(-\infty)$-tuple.  This means that if $\vec{u}\in\{{\rm Standard\,\,forms\,\,for\,\,}\vec{s}\cdot\vec{t}\,\,|\,\vec{s}\in 0',\,\,\vec{t}\in 0'\}$, then $\vec{u}\in 0'$.  Furthermore, if $\vec{r}\in 0'$, then we can write $\vec{r}$ as 
$$\begin{array}{rcl}
\vec{r}&=&(s_x^*,s_y^*,...\,,s_z^*,...) \cr
&=&(s_x^*,...\,,s_u^*,...\,, s_v^*,...\,,s_z^*,...) \cr
&=&(s_x^*,...\,,s_u^*,...)\cdot(s_v^*,...\,,s_z^*,...). \cr
\end{array}$$
This means that
$$\vec{r}=(s_x^*,...\,,s_u^*,...)\cdot(s_v^*,...\,,s_z^*,...)$$
where tuples $(s_x^*,...\,,s_u^*,...)\in 0'$ and $(s_v^*,...\,,s_z^*,...)\in 0'$.  Therefore, we see that $\vec{r}\in\{{\rm Standard\,\,forms\,\,for\,\,}\vec{s}\cdot\vec{t}\,\,|\,\vec{s}\in 0',\,\,\vec{t}\in 0'\}$.  Putting this all together, we see
$$\begin{array}{rcl}
0'\cdot\,0'&=&\{{\rm Standard\,\,forms\,\,for\,\,}\vec{s}+\vec{t}\,\,|\,\vec{s}\in 0',\,\,\vec{t}\in 0'\} \cr
&=&0'. \cr
\end{array}$$
We conclude that $0'\cdot\,0'=0'$.
\newline \indent If however $[n,x]\in[\Bbb Z,\Bbb S_0]$, then we see that
$$[n,x]\cdot \,0'=\{{\rm Standard\,\,forms\,\,for\,\,}\vec{s}\cdot\vec{t}\,\,|\,\vec{s}\in [n,x],\,\,\vec{t}\in 0'\}.$$
In Section 2, we saw that $\vec{s}\cdot\vec{t}$ will also be a $(-\infty)$-tuple.  This means that if $\vec{u}\in\{{\rm Standard\,\,forms\,\,for\,\,}\vec{s}\cdot\vec{t}\,\,|\,\vec{s}\in [n,x],\,\,\vec{t}\in 0'\}$, then $\vec{u}\in 0'$.  Furthermore, if $\vec{r}\in 0'$, then we can write $\vec{r}$ as
$$\begin{array}{rcl}
\vec{r}&=&(s_x^*,s_1^*,s_1^*,...\,,s_1^*,s_y^*,...\,,s_z^*,...) \cr
&=&(s_x^*,s_1^*,...\,,s_1^*)\cdot(s_y^*,...\,,s_z^*,...) \cr
\end{array}$$
where tuples $(s_x^*,s_1^*,...\,,s_1^*)\in[n,x]$ and $(s_y^*,...\,,s_z^*,...)\in 0'$.  Therefore, we see that $\vec{r}\in\{{\rm Standard\,\,forms\,\,for\,\,}\vec{s}\cdot\vec{t}\,\,|\,\vec{s}\in [n,x],\,\,\vec{t}\in 0'\}$.  Putting this all together, we see
$$\begin{array}{rcl}
[n,x]\cdot 0'&=&\{{\rm Standard\,\,forms\,\,for\,\,}\vec{s}\cdot\vec{t}\,\,|\,\vec{s}\in [n,x],\,\,\vec{t}\in 0'\} \cr
&=&0'. \cr
\end{array}$$
Therefore $[n,x]\cdot 0'=0'$.  Similarly, we can show that $0'\cdot [n,x]=0'$.  We conclude that $[n,x]\
\cdot 0'=0'\cdot [n,x]=0'$ for every $[n,x]\in[\Bbb Z,\Bbb S_0]$.  We also conclude that absolute zero multiplies with elements in $[\Bbb Z,\Bbb S_0]^*$ just like the additive identity in a field multiplies with other elements in the field.

\section*{Finding a Value for Absolute Zero}
Recall from Section 2 that the set $0'$ only has level $-\infty$.  But we can find a value for $0'$ as well if we express absolute zero much like we express elements in $[\Bbb Z,\Bbb S_0]$, as $0'=[-\infty,x]$ with value $x\in\Bbb S_0\setminus\{0\}$.  We know that $0'+\,0'=0'$.  And if we assume the formula we derived for addition of elements in $[\Bbb Z,\Bbb S_0]$ applies for absolute zero, then we see that
$$\begin{array}{rcl}
{\rm [}-\infty,x{\rm ]}&=&0' \cr
&=&0'+\,0' \cr
&=&{\rm [}-\infty,x{\rm ]}+{\rm [}-\infty,x{\rm ]} \cr
&=&{\rm [}-\infty,x+x{\rm ]} \cr
&=&{\rm [}-\infty,2x{\rm ]}. \cr
\end{array}$$
Therefore $[-\infty,x]=[-\infty,2x]$, which implies that $x=2x$.  Now there is no non-zero $x\in\Bbb S_0$ such that $x=2x$ as $(\Bbb S_0,+,\cdot)$ is a field.  But if we relax the requirement that $x$ be non-zero and allow $x=0$, then we have a solution.  So we say that 
$$0'=[-\infty,0]$$
with value $0\in\Bbb S_0$.  Note that $0'$ is the only element in $[\Bbb Z,\Bbb S_0]^*$ with a value of $0$.  

\section*{Complete S-Extensions of Fields}
We have seen how to add and multiply together any two elements in $[\Bbb Z,\Bbb S_0]^*$.  In particular, we know how to add elements in $[\Bbb Z,\Bbb S_0]^*$ if they have different levels.  We know that $0'+\,0'=0'$ and if $[n,x],[m,y]\in[\Bbb Z,\Bbb S_0]$, then we also know how to evaluate $[n,x]+[m,y]$ whenever $m=n$, provided $x+y\not=0$.  But now we can use absolute zero to define addition of elements $[n,x],[m,y]\in[\Bbb Z,\Bbb S_0]^*$ where $n=m$ and $x+y=0$ as  
$$[n,x]+[m,y]=0'.$$
With that, we have created the structure $([\Bbb Z,\Bbb S_0]^*,+,\cdot)$ which we will call the {\bf Complete S-Extension Of The Field $(\Bbb S_0,+,\cdot)$} with the following addition and multiplication operations defined for all possible $[n,x],[m,y]\in[\Bbb Z,\Bbb S_0]^*$:
$$\begin{array}{rcl}
[n,x]+[m,y]&=&\left\{
				\begin{array}{cl}
						{\rm [}n,x+y{\rm ]} & \quad {\rm if}\,\,n=m,\,\,x+y\not=0 \\
						0' & \quad {\rm if}\,\,n=m,\,\,x+y=0 \\
						{\rm [}n,x{\rm ]} & \quad {\rm if}\,\,n>m \\
						{\rm [}m,y{\rm ]} & \quad {\rm if}\,\,n<m \\
        \end{array} \right. \cr
&& \cr
{\rm [}n,x{\rm ]}\cdot {\rm [}m,y{\rm ]}&=&{\rm [}n+m, x\cdot y{\rm ]}. \cr				
\end{array}$$
Assuming $-\infty+(-\infty)=-\infty$ and $-\infty+m=m+(-\infty)=-\infty$ for every $m\in\Bbb Z$, the element $0'=[-\infty,0]$ can be considered when using these operations and will yield $0'+[n,x]=[n,x]$ and $0'\cdot [n,x]=0'$ for every $[n,x]\in[\Bbb Z,\Bbb S_0]$ while $0'+\,0'=0'$ and $0'\cdot \,0'=0'$, as was shown earlier.  
\newline \indent Subtraction is defined in accordance with subtraction in full s-extensions so that $[n,x]-[m,y]=[n,x]+(-[m,y])$ for every $[n,x],[m,y]\in[\Bbb Z, \Bbb R]^*$.  And since $-[m,y]=[m,-y]$, this means that
$$\begin{array}{rcl}
[n,x]-[m,y]&=&[n,x]+(-[m,y]) \cr
           &=&[n,x]+[m,-y] \cr
					 &=&\left\{
						\begin{array}{cl}
								{\rm [}n,x-y{\rm ]} & \quad {\rm if}\,\,n=m,\,\,x-y\not=0 \\
								0' & \quad {\rm if}\,\,n=m,\,\,x-y=0 \\
								{\rm [}n,x{\rm ]} & \quad {\rm if}\,\,n>m \\
								{\rm [}m,-y{\rm ]} & \quad {\rm if}\,\,n<m \\
						\end{array} \right. \cr
\end{array}$$
Most importantly, for every $[n,x]\in[\Bbb Z,\Bbb S_0]^*$, we see that
$$[n,x]-[n,x]=[n,x]+[n,-x]=0'.$$
This includes the case where $[n,x]=[-\infty,0]=0'$.  To differentiate this subtraction from the subtraction in $([\Bbb Z,\Bbb S_0],+,\cdot)$ which cannot be defined due to its indeterminate nature, we call this subtraction {\bf Absolute Subtraction}.
\newline \indent Recall that division of elements $[n,x]$ and $[m,y]$ in a full s-extension is defined as 
$$\frac{[n,x]}{[m,y]}=q\,\,\,{\rm such\,\,that}\,\,\,[m,y]\cdot q=[n,x].$$
We can use this same definition to define division in a complete s-extension.  In this way, we know that $\frac{[n,x]}{[m,y]}=[n-m,\frac{x}{y}]$ for all $[n,x],[m,y]\in[\Bbb Z,\Bbb S_0]$.  But if $[m,y]\in[\Bbb Z,\Bbb S_0]^*$ is not absolute zero and we let $\frac{0'}{[m,y]}=[p,z]$ where $[p,z]\in[\Bbb Z,\Bbb S_0]^*$, then $[m,y]\cdot[p,z]=0'$.  Therefore $[m+p,y\cdot z]=[-\infty,0]$ implying that $m+p=-\infty$ and $y\cdot z=0$.  However $[m,y]\not=[-\infty,0]$, which means $m\in\Bbb Z$ and $y\in\Bbb S_0\setminus\{0\}$.  So $m+p=-\infty$ only if $p=-\infty$ as $m\in\Bbb Z$ while $y\cdot z=0$ only if $z=0$ as $y\not=0$.  This means $\frac{0'}{[m,y]}=0'$ as  
$$\frac{0'}{[m,y]}=\frac{[-\infty,0]}{[m,y]}=[p,z]=[-\infty,0]=0'.$$
Furthermore, notice that if $[n,x]=[-\infty,0]$, then
$$\frac{[n,x]}{[m,y]}=\frac{[-\infty,0]}{[m,y]}=[-\infty,0]=[-\infty-m,\frac{0}{y}]=[n-m,\frac{x}{y}].$$
So we find that $\frac{[n,x]}{[m,y]}=[n-m,\frac{x}{y}]$ in the case where $[n,x]=[-\infty,0]$.  This means division satisfies $\frac{[n,x]}{[m,y]}=[n-m,\frac{x}{y}]$ for every $[n,x]\in[\Bbb Z,\Bbb S_0]^*$ and every $[m,y]\in[\Bbb Z,\Bbb S_0]^*$ that is not absolute zero.  
\newline \indent We can complete our understanding of division by considering the case of $\frac{0'}{0'}$ and the case of $\frac{[n,x]}{0'}$.  If we let $\frac{0'}{0'}=(p,z)$, then $0'=0'\cdot [p,z]$.  Since every $[p,z]\in[\Bbb Z,\Bbb S_0]^*$ satisfies this condition, we say that $\frac{0'}{0'}$ is undefined as it s indeterminate.  At the same time, if $[n,x]\in[\Bbb Z,\Bbb S_0]$ is not absolute zero, then we see that $\frac{[n,x]}{0'}$ must satisfy $\frac{[n,x]}{0'}=[p,z]$ such that $0'\cdot[p,z]=[n,x]$.  However $0'\cdot[p,z]=0'\not=[n,x]$, implying that there is no element $[p,z]\in[\Bbb Z,\Bbb S_0]^*$ that satisfies this condition.  We say that $\frac{[n,x]}{0'}$ is undefined as a solution in $[\Bbb Z,\Bbb S_0]^*$ does not exist.  Putting this together, we see that division is defined for every $[n,x],[m,y]\in[\Bbb Z,\Bbb S_0]^*$ where $[m,y]$ is not absolute zero and satisfies
$$\frac{[n,x]}{[m,y]}=[n-m,\frac{x}{y}]$$
while division of the form $\frac{[n,x]}{0'}$ for any element $[n,x]\in[\Bbb Z,\Bbb S_0]^*$ is undefined.
\newline \indent Before moving on, it is worth pointing out that in $([\Bbb Z,\Bbb S_0]^*,+,\cdot)$, the arithmetic operations of addition, subtraction, multiplication, and division are well defined for every element in $[\Bbb Z,\Bbb S_0]^*$ with the only exception being division by absolute zero.  It's also worth proving that the complete s-extension of a field is unique to the field it extends.  

\begin{proposition}
Let $(\Bbb S_1,+_1,\cdot_1)$ and $(\Bbb S_2,+_2,\cdot_2)$ be ascended fields with the scalar fields $(\Bbb S_0^1,+_1,\cdot_1)$ and $(\Bbb S_0^2,+_2,\cdot_2)$, respectively.  If $(\Bbb S_0^1,+_1,\cdot_1)=(\Bbb S_0^2,+_2,\cdot_2)$, then
$$([\Bbb Z,\Bbb S_0]^*_1,+_1,\cdot_1)=([\Bbb Z,\Bbb S_0]^*_1,+_1,\cdot_1).$$

\end{proposition}

\begin{proof}
We have previously shown that $([\Bbb Z,\Bbb S_0]_1,+_1,\cdot_1)=([\Bbb Z,\Bbb S_0]_2,+_2,\cdot_2)$.  Therefore $[\Bbb Z,\Bbb S_0]_1=[\Bbb Z,\Bbb S_0]_2$.  Furthermore $\Bbb S_0^1=\Bbb S_0^2$, which means $0'_1=0'_2$ as 
$$0'_1=[-\infty,0_1]=[\infty,0_2]=0'_2.$$
This means $[\Bbb Z,\Bbb S_0]_1\cup\{0'_1\}=[\Bbb Z,\Bbb S_0]_2\cup\{0'_2\}$, implying $[\Bbb Z,\Bbb S_0]^*_1=[\Bbb Z,\Bbb S_0]^*_2$.  
\newline \indent Since $([\Bbb Z,\Bbb S_0]_1,+_1,\cdot_1)=([\Bbb Z,\Bbb S_0]_2,+_2,\cdot_2)$, $[n,x]+_1[m,y]=[n,x]+_2[m,y]$ and $[n,x]\cdot_1[m,y]=[n,x]\cdot_2[m,y]$ for every $[n,x],[m,y]\in[\Bbb Z,\Bbb S_0]_1=[\Bbb Z,\Bbb S_0]_2$.  Now $0'_1+_1 0'_1=0'_2+_2 0'_2$ as 
$$0'_1+_1 0'_1=0'_1=0'_2=0'_2+_2 0'_2$$
while $0'_1\cdot_1 0'_1=0'_2\cdot_2 0'_2$ as 
$$0'_1\cdot_1 0'_1=0'_1=0'_2=0'_2\cdot_2 0'_2.$$
We also see that $0'_1+_1 [n,x]=0'_2+_2 [n,x]$ for every $[n,x]\in[\Bbb Z,\Bbb S_0]_1=[\Bbb Z,\Bbb S_0]_2$ as 
$$0'_1+_1 [n,x]=[n,x]=0'_2+_2 [n,x]$$
while $0'_1\cdot_1 [n,x]=0'_2\cdot_2 [n,x]$ for every $[n,x]\in[\Bbb Z,\Bbb S_0]_1=[\Bbb Z,\Bbb S_0]_2$ as 
$$0'_1\cdot_1 [n,x]=0'_1=0'_2=0'_2\cdot_2 [n,x].$$
Putting this all together, we see that $[n,x]+_1[m,y]=[n,x]+_2[m,y]$ and $[n,x]\cdot_1[m,y]=[n,x]\cdot_2[m,y]$ for every $[n,x],[m,y]\in[\Bbb Z,\Bbb S_0]^*_1=[\Bbb Z,\Bbb S_0]^*_2$.  Therefore we can conclude that $([\Bbb Z,\Bbb S_0]^*_1,+_1,\cdot_1)=([\Bbb Z,\Bbb S_0]^*_1,+_1,\cdot_1).$ $\square$ 
\end{proof} 

\begin{corollary}
Let $(F,+_F,\cdot_F)$ be any field, then there exists a unique complete s-extension of $(F\times\{0\},+,\cdot)$ where $(F\times\{0\},+,\cdot)$ is the scalar field of an ascended field $(F\times F,+,\cdot)$ and is isomorphic to $(F,+_F,\cdot_F)$.
\end{corollary}

\begin{proof}
We have previously shown that for every field $(F,+_F,\cdot_F)$ there exists a unique full s-extension $([\Bbb Z,F\times\{0\}],+,\cdot)$ where $(F\times\{0\},+,\cdot)$ is the scalar field of an ascended field $(F\times F,+,\cdot)$ and is isomorphic to $(F,+_F,\cdot_F)$.  By considering the absolute zero of $([\Bbb Z,F\times\{0\}],+,\cdot)$, we can create the complete s-extension $([\Bbb Z,F\times\{0\}]^*,+,\cdot)$.  And by the previous proposition, this complete s-extension must be the same as any other complete s-extension built from an ascended field $(F\times F,+,\cdot)$ containing the scalar field $(F\times\{0\},+,\cdot)$.  Therefore $([\Bbb Z,F\times\{0\}]^*,+,\cdot)$ is the unique complete s-extension of $(F\times\{0\},+,\cdot)$. $\square$
\end{proof}
In light of this corollary, $([\Bbb Z,F],+,\cdot)$, the complete s-extension of the field $(F,+_F,\cdot_F)$, and $([\Bbb Z,F\times\{0\}],+,\cdot)$, the complete s-extension of the scalar field $(F\times\{0\},+,\cdot)$ of an ascended field $(F\times F,+,\cdot)$, are often used interchangeably.  
\newline \indent This proposition can also be generalized to two ascended fields $(\Bbb S_1,+_1,\cdot_1)$ and $(\Bbb S_2,+_2,\cdot_2)$ containing the fields $(\Bbb S_0^1,+_1,\cdot_1)$ and $(\Bbb S_0^2,+_2,\cdot_2)$ which are isomorphic instead of equal.  Since $(\Bbb S_0^1,+_1,\cdot_1)\cong(\Bbb S_0^2,+_2,\cdot_2)$, there exists an isomorphism $f:\Bbb S_0^1\to\Bbb S_0^2$.  We then define the function $\phi:[\Bbb Z, \Bbb S_0^1]^*\to[\Bbb Z, \Bbb S_0^2]^*$ such that $\phi(0'_1)=0'_2$ and $\phi(s)=[n,f(x)]$ for all $s=[n,x]\in[\Bbb Z, \Bbb S_0^1]^*$ that are not absolute zero.  Now this function $\phi$ is an isomorphism.  So if $(\Bbb S_1,+_1,\cdot_1)$ and $(\Bbb S_2,+_2,\cdot_2)$ are two ascended fields containing scalar fields $(\Bbb S_0^1,+_1,\cdot_1)$ and $(\Bbb S_0^2,+_2,\cdot_2)$, respectively, where $(\Bbb S_1,+_1,\cdot_1)\cong(\Bbb S_2,+_2,\cdot_2)$, then $([\Bbb Z, \Bbb S_0^1]^*,+_1,\cdot_1)\cong([\Bbb Z, \Bbb S_0^2]^*,+_2,\cdot_2)$. 
\newline \indent Finally, we can show that there is an isomorphism from $\Bbb S_0$ to $[0,\Bbb S_0]^*$ where
$$[0,\Bbb S_0]^*=\big\{\,[0,x]\,\,|\,x\in\Bbb S_0\setminus\{0\}\big\}\,\,\cup\,\,\big\{0'\big\}$$
by defining the function $\phi:\Bbb S_0\rightarrow [0,\Bbb S_0]^*$ where $\phi(0)=0'$ and $\phi(x)=[0,x]$ for all $x\not=0$.  It is easy to show that this function is a bijection.  And if $x,y\in\Bbb S_0$ are non-zero, then we can see that 
$$\phi(x)+\phi(y)=[0,x]+[0,y]=[0,x+y]=\phi(x+y)$$
$$\phi(x)\cdot \phi(y)=[0,x]\cdot[0,y]=[0+0,x\cdot y]=[0,x\cdot y]=\phi(x\cdot y).$$
At the same time, we see that
$$\phi(0)+\phi(y)=[-\infty,0]+[0,y]=[0,y]=\phi(y)=\phi(0+y)$$
$$\phi(0)\cdot\phi(y)=[-\infty,0]\cdot [0,y]=[-\infty,0]=\phi(0)=\phi(0\cdot y)$$
$$\phi(0)+\phi(0)=[-\infty,0]+[-\infty,0]=[-\infty,0]=\phi(0)=\phi(0+0)$$
$$\phi(0)\cdot\phi(0)=[-\infty,0]\cdot [-\infty,0]=[-\infty,0]=\phi(0)=\phi(0\cdot 0).$$
Therefore $\phi$ is an isomorphism, allowing us to recover the field $\Bbb S_0$ from the complete s-extension of the field by considering $[0,\Bbb S_0]^*$.  It is worth noting that in this isomorphism, the field element $0$ is mapped to absolute zero $0'$ in the complete s-extension of the field.

\section*{Order in the Complete S-Extension of the Reals}
Let $([\Bbb Z,\Bbb R],+,\cdot)$ be the complete s-extension of the real numbers $(\Bbb R,+,\cdot)$.  Using the set of positive real numbers $\Bbb P(\Bbb R)$, we can define the set $\Bbb P$ as
$$\Bbb P=\{[n,x]\,|\,x\in\Bbb P(\Bbb R)\}.$$
Note that $[-\infty,0]\notin\Bbb P$ as $0\notin\Bbb P(\Bbb R)$.  Now, if $[n,x],[m,y]\in\Bbb P$, then $x,y\in\Bbb P(\Bbb R)$.  We then see that 
$$[n,x]+[m,y]=\left\{
				\begin{array}{ll}
						{\rm [}n,x+y{\rm ]} & \quad {\rm if}\,\,n=m \\
						{\rm [}n,x{\rm ]} & \quad {\rm if}\,\,n>m \\
						{\rm [}m,y{\rm ]} & \quad {\rm if}\,\,n<m \\
        \end{array} \right.$$
where $x+y\in\Bbb P(\Bbb R)$.  We also see that $[n,x]\cdot [m,y]=[n+m,x\cdot y]$ where $x\cdot y\in\Bbb P(\Bbb R)$.  Therefore $[n,x]+[m,y]\in\Bbb P$ and $[n,x]\cdot [m,y]\in\Bbb P$.  Additionally, if $[n,x]\in[\Bbb Z,\Bbb R]^*$, then we see that either $[n,x]\in\Bbb P$, $-[n,x]=[n,-x]\in\Bbb P$, or $[n,x]=[-\infty,0]$ as either $x\in\Bbb P(\Bbb R)$, $-x\in\Bbb P(\Bbb R)$, or $x=0$.  Hence the complete s-extension of the reals satisfies the property of trichotomy.  
\newline \indent We can then setup a natural ordering for elements $[n,x],[m,y]\in[\Bbb Z,\Bbb R]^*$ as
$$[n,x]<[m,y]\,\,{\rm if\,\,and\,\,only\,\,if}\,\,[m,y]-[n,x]\in\Bbb P.$$
Using this ordering, we find that for the complete s-extension $([\Bbb Z,\Bbb R]^*,+,\cdot)$, we have
$$...\,<-[n,x]<-[0,z]<-[-m,y]<0'<[-m,y]<[0,z]<[n,x]<...$$
where $[n,x],[0,z],[-m,y]\in\Bbb P$ with $n,-m\in\Bbb Z$ where $n$ is positive and $-m$ is negative.  This can also be written as
$$...\,<[n,-x]<[0,-z]<[-m,-y]<0'<[-m,y]<[0,z]<[n,x]<...$$
where $x,y,z\in\Bbb P(\Bbb R)$.  Note that any element with a positive value is a positive number and any element with a negative value is a negative number while absolute zero, $0'$, separates the positive numbers and the negative numbers.  
\newline \indent Most importantly though, with this ordering, we can now consider any $[0,x]\in[\Bbb Z,\Bbb R]$ to be a {\bf Finite Real Number}, any element $[n,x]\in[\Bbb Z,\Bbb R]$ where $n>0$ to be an {\bf Infinite Real Number of Level n}, and any $[-n,x]\in[\Bbb Z,\Bbb R]$ where $-n<0$ to be an {\bf Infinitesimal Real Number of Level n}.  Note that infinite real numbers are either greater than or less than every possible finite real number while infinitesimal real numbers are not absolute zero but are less than every possible positive finite real number while simultaneously being greater than every possible negative finite real number.  With this, we have introduced finite, infinite, and infinitesimal real numbers in a way that the collection $[0,\Bbb R]^*$ of finite real numbers and absolute zero corresponds to the real numbers $\Bbb R$.

\section*{Finite, Infinite, and Infinitesimal Real Numbers}
Notice that for any non-zero, real number $x\in\Bbb R$, we can let $x$ be the value of a number in $[\Bbb Z,\Bbb R]$.  In a sense, we can then choose the numerical form that this value $x$ assumes when expressed as a number in $[\Bbb Z,\Bbb R]$ by assigning it a level.  If the level is $0$, then we can think of $x$ as taking on its finite form $[0,x]$.  If the level is $1$, then we can think of $x$ as taking on its first infinite form $[1,x]$, which is commonly written as $[x]$ from the set $[\Bbb R]$ which we talked about in Section 3.  And if the level is $-1$, then we can think of $x$ as taking on its first infinitesimal form $[-1,x]$, which is commonly written as $\epsilon_x$ from the set $\epsilon_{\Bbb R}$ where
$$\epsilon_{\Bbb R}=\{\epsilon_x\,|\,x\in\Bbb R\setminus\{0\}\}$$
such that $\epsilon_x=[-1,x]\in[\Bbb Z,\Bbb R]$ for each $x\in\Bbb R\setminus\{0\}$.  
\newline \indent In fact, every $[n,x]\in[\Bbb Z,\Bbb R]$ with $n>0$ can be thought of as the nth infinite form of $x$ while every $[-n,x]$ with $-n<0$ can be thought of a the nth infinitesimal form of $x$.  In this way, we justify the name finite, infinite, and infinitesimal real numbers as each $[n,x]\in[\Bbb Z,\Bbb R]$ with a fixed $x\in\Bbb R\setminus\{0\}$ can be thought of as a either a finite, infinite, or infinitesimal form for the non-zero, real number value $x$.  We can also justify calling any number $[0,x]\in[\Bbb Z,\Bbb R]$ the {\bf base form for the value $x$} as $[0,x]$ corresponds to the real number $x$ in its common or base understanding that all mathematicians are accustomed to. 
\newline \indent It is worth noting that multiplying an infinite real number $[n,x]\in[\Bbb Z,\Bbb R]$ with level $n$ and value $x$ by an infinitesimal real number $[-n,y]\in[\Bbb Z,\Bbb R]$ with level $-n$ and value $y$ yields
$$[n,x]\cdot[-n,y]=[n+(-n),x\cdot y]=[0,x\cdot y].$$
Hence $[n,x]\cdot[-n,y]=[0,x\cdot y]$ which is a finite real number corresponding to the real number $x\cdot y$.  Therefore, multiplying an infinite real number of level $n$ by an infinitesimal real number of level $-n$ yields a finite real number.  
\newline \indent Furthermore, if we take a finite real number $[0,x]\in[\Bbb Z,\Bbb R]$ with value $x$ and divide it by an infinite real number $[n,y]\in[\Bbb Z,\Bbb R]$ with level $n$ and value $y$, then the result is
$$\frac{[0,x]}{[n,y]}=[0-n,\frac{x}{y}]=[-n,\frac{x}{y}]$$
which is an infinitesimal real number with level $-n$ and value $\frac{x}{y}$.  If a finite real number $[0,x]\in[\Bbb Z,\Bbb R]$ is divided by an infinitesimal real number $[-n,y]\in[\Bbb Z,\Bbb R]$ with level $-n$ and value $y$, then the result is
$$\frac{[0,x]}{[-n,y]}=[0-(-n),\frac{x}{y}]=[n,\frac{x}{y}]$$
which is an infinite real number with level $n$ and value $\frac{x}{y}$.  Therefore dividing a finite real number by an infinite real number results in an infinitesimal real number while dividing a finite real number by an infinitesimal real number results in an infinite real number.

\section{Future Works and Applications}
With the development of infinite and infinitesimal real numbers, it is proposed that the following function $f:[\Bbb Z,\Bbb R]^*\times[\Bbb Z,\Bbb R]^*\longrightarrow[\Bbb Z,\Bbb R]^*$ where 
$$f([n,x],[m,y])=\big|\,[n,x]-[m,y]\,\big|$$
be used to measure distance in $[\Bbb Z,\Bbb R]^*$.  The author has shown this function to satisfy most of the properties of a metric, with the proof showing that the triangle inequality holds requiring final edits.  
\newline \indent Conditioned on $f$ being a metric, two possible applications are being considered.  The first involves understanding the output of real functions at infinitesimal distances around possible real inputs in connection with the results of division by zero in the simple s-extension of the reals applied to the first infinite forms of the real numbers for functions with division by zero singularities.  
The second seeks a geometric interpretation for the complete s-extension of the reals analogous to the real number line interpretation of the real numbers and considers line segments with infinitesimal lengths, which the author calls {\bf pixels}.  Ongoing investigations into these applications are currently underway. 
\newline \indent Finally, the author has shown that the complete s-extension of any field satisfies all of the requirements to be a field except additive associativity and will release a manuscript with the accompanying proof as a supplement to this manuscript. 


\end{document}